\documentclass[12pt]{amsart}
\usepackage{dsfont,hyperref}
\usepackage[all]{xy}
\usepackage{amssymb,hyperref}
\SelectTips{cm}{}
\allowdisplaybreaks
\numberwithin{equation}{section}
\usepackage{color}

\usepackage{graphicx}
\usepackage{amsmath,amsthm, amscd, amssymb, amsfonts}
\usepackage{enumitem}
\usepackage[all]{xy}
\usepackage{tikz-cd}
\usepackage{color}
\usepackage{xcolor}
\usepackage{nicefrac}
\usepackage[bottom]{footmisc}
\usepackage{relsize}
\usepackage{latexsym}
\usepackage{hyperref}
\usepackage{comment}
\usepackage{mathrsfs}
\usepackage[enableskew]{youngtab}
\usepackage[all]{xy}
\usepackage{empheq}
\usepackage{empheq}
\usepackage{tikz}
\usetikzlibrary{matrix}
\usepackage{combelow}

\newtheorem{theorem}{Theorem}[section]
\newtheorem{lemma}[theorem]{Lemma}
\newtheorem{proposition}[theorem]{Proposition}
\newtheorem{corollary}[theorem]{Corollary}
\newtheorem{conjecture}[theorem]{Conjecture}

\newtheorem{thm}[theorem]{Theorem}
\newtheorem{rem}[theorem]{Remark}
\newtheorem{defn}[theorem]{Definition}
\newtheorem{lem}[theorem]{Lemma}
\theoremstyle{plain}

\newcommand\numeq[1]%
{\stackrel{\scriptscriptstyle(\mkern-1.5mu#1\mkern-1.5mu)}{=}}

\newcommand{\id}{\text{id}}

\newcommand{\rank}{\text{rank\,}}

\newcommand{\Rep}{\text{Rep}}

\newcommand{\eps}{\varepsilon}

\newcommand{\ad}{{\text{ad}}}

\newcommand{\lk}{\ar@{-}}

\DeclareMathOperator{\ev}{ev}

\newcommand{\ncm}{\newcommand}

\ncm{\blue}{\textcolor[rgb]{.00, .00, 1.00}}
\ncm{\red}{\textcolor[rgb]{1.00, .00, .00}}
\ncm{\green}{\textcolor[rgb]{.50, 0.20, .90}}

\ncm{\np}{\newpage}
\ncm{\ebl}{\end{thebibliography}}
\ncm{\bbl}{\begin{thebibliography}}
\ncm{\chd}{_{ _{\ch}}}
\ncm{\ald}{_{ _{\al}}}
\newcommand{\blam}{\Lam}
\ncm{\cP}{\mathcal{P}}
\ncm{\ei}{e_i}
\ncm{\eij}{e_{i,\;j}}
\ncm{\bt}{\begin{thm}}
\ncm{\bdef}{\begin{defn}}
\ncm{\edf}{\end{defn}}
\ncm{\et}{\end{thm}}
\ncm{\bc}{\begin{corollary}}
\ncm{\bl}{\begin{lem}}
\ncm{\el}{\end{lem}}
\ncm{\bpf}{\begin{proof}}
\ncm{\epf}{\end{proof}}
\ncm{\ec}{\end{corollary}}
\ncm{\ord}{\mtr{ord}}
\ncm{\er}{\end{rem}}
\ncm{\br}{\begin{rem}}
\ncm{\bn}{\begin}
\ncm{\bp}{\begin{proposition}}
\ncm{\ep}{\end{proposition}}
\ncm{\bd}{\begin{document}}
\ncm{\ed}{\end{document}}
\ncm{\beq}{\begin{equation}}
\ncm{\beqn}{\begin{equation*}}
\ncm{\eeq}{\end{equation}}
\ncm{\eeqn}{\end{equation*}}
\ncm{\bea}{\begin{eqnarray}}
\ncm{\eea}{\end{eqnarray}}
\ncm{\beanon}{\begin{eqnarray*}}
\ncm{\eeanon}{\end{eqnarray*}}\ncm{\ek}{\eps|_K}\ncm{\diez}{\#}
\ncm{\bwt}{\bowtie}
\ncm{\cC}{\mtc{C}}\ncm{\cc}{\mtc{C}}
\ncm{\cX}{\mtc{X}}
\ncm{\wt}{\widetilde}
\ncm{\sg}{\sigma}
\ncm{\X}{\mathcal{X}}
\ncm{\cA}{\mathcal{A}}
\ncm{\HKer}{\mtr{HKer}}
\ncm{\LKER}{\mtr{LKer}}
\ncm{\aad}{\mtr{ad}}
\newcommand{\mbf}{\mathbb F}
\ncm{\Dr}{\mtr{D}}
\ncm{\cD}{{\mathcal{D}}}\ncm{\cd}{{\mathcal{D}}}\ncm{\ce}{{\mathcal{E}}}
\ncm{\G}{\mathcal{G}}
\ncm{\Dc}{\mtc{D}}
\ncm{\fp}{\mathsf{FPdim}}
\ncm{\Vc}{\mtr{Vec}}
\ncm{\cK}{\mtc{K}}
\ncm{\cM}{\mtc{M}}
\ncm{\cE}{\mtc{E}}
\ncm{\cS}{\mtc{S}}
\ncm{\cop}{\mtr{cop}}
\ncm{\chr}{character }\ncm{\ck}{\mtc{K}}
\ncm{\bw}{\bwt}
\ncm{\hker}{\mtr{HKer}}
\ncm{\bx}{\boxtimes}

\ncm{\bne}{\begin{enumerate}}
\ncm{\ene}{\end{enumerate}}
\ncm{\lker}{\mtr{LKer}}
\ncm{\md}{\medbreak}
\ncm{\rep}{\Rep}
\ncm{\ind}{\mtr{ind}}
\ncm{\mdn}{\md\noindent}
\ncm{\dd}{$}
\ncm{\up}{^}

\newcommand\bS{\mathbf{S}}
\newcommand\RR{{\mathbb R}} 
\newcommand\FF{{\mathbb F}}
\newcommand{\fpcd}{{\fp(\cd)}}
\newcommand\htimes{\,\widehat\otimes\,}
\newcommand\pib{\overline{\pi}}  
\newcommand\groth{{G_0(A)}}
\newcommand\xx{\mathbf{x}}
\newcommand\sss{\mathbf{s}}
\newcommand\ppp{\mathbf{p}}
\newcommand\nnn{\mathbf{n}}
\newcommand\GG{\mathbf{G}}
\newcommand\NN{\mathbb{N}}
\newcommand{\ccb}{{\mathcal B}}
\newcommand{\bdelta}{\Delta}
\newcommand{\bann}{\mtr{Ann}}
\newcommand{\coh}{{\rm H}}
\newcommand{\exx}{\mtr{End}}
\newcommand{\HH}{{\rm HH}}
\newcommand{\ra}{\rightarrow}
\newcommand{\codim}{{\rm codim}}
\newcommand{\Wedge}{\textstyle\bigwedge}
\newcommand{\Z}{\mathbb Z}
\newcommand{\ot}{\otimes}
\newcommand{\co}{\mathcal O}
\newcommand{\xra}{\xrightarrow}
\newcommand{\mtc}{\mathcal}
\newcommand{\cs}{\mtc S}
\newcommand{\onh}{On the other hand}
\newcommand{\lam}{\lambda}
\newcommand{\lb}{\label}
\newcommand{\Lam}{\Lambda}
\newcommand{\lbd}{\Lambda}
\newcommand{\ann}{\mtr{Ann}}
\newcommand{\sig}{\sigma}
\newcommand{\al}{\alpha}
\newcommand{\en}{\end}
\newcommand{\wte}{{\widehta{E}}}
\newcommand{\rest}{\mtr{Res}}
\newcommand{\sub}{\subsection}
\newcommand{\mc}{\mathcal}
\newcommand{\teta}{\theta}
\newcommand{\nhs}{normal Hopf subalgebra}
\newcommand{\ul}{\underline}
\newcommand{\mr}{\mathrm}
\newcommand{\rh}{\rightharpoonup}
\newcommand{\lh}{\leftharpoonup}
\newcommand{\whb}{{\widehat{(H, \mtc B)}}}
\newcommand{\hb}{{(H, \mtc B)}}
\newcommand{\fjss}{F^j_{ss}}
\newcommand{\csujss}{\mtr{C}^j_{ss}}
\newcommand{\csujts}{\mtr{C}^j_{ts}}
\newcommand{\sprime}{{s'}}
\newcommand{\pt}{{\mtr{pt}}}
\newcommand{\pif}{\pi_f}
\newcommand{\B}{\mtc{B}}
\newcommand{\fjts}{F^j_{ts}}
\newcommand{\iotal}{\iota_{\mthl}}
\newtheorem{exmp}{Example}[section]
\newcommand{\hmcd}{{\hm_\cd}}
\newcommand{\csuj}{\csu^j}
\newcommand{\mujst}{\mu^j_{st}}
\newcommand{\mukuv}{\mu^k_{uv}}
\newcommand{\mujpsptp}{\mu^{j'}_{s't'}}
\newcommand{\cjst}{\csu^j_{st}}
\newcommand{\cjts}{\csu^j_{ts}}
\newcommand{\cjpsptp}{\csu^{j'}_{s't'}}
\newcommand{\cjptpsp}{\csu^{j'}_{t's'}}
\newcommand{\ncjst}{\frac{\cjst}{\dim(\cc_j)}}
\newcommand{\ncjts}{\frac{\cjts}{\dim(\cc_j)}}
\newcommand{\ncjpsptp}{\frac{\cjpsptp}{\dim(\cc_{j'})}}
\newcommand{\ncjptpsp}{\frac{\cjptpsp}{\dim(\cc_{j'})}}
\newcommand{\lu}{\mtc L(\unu)}
\newcommand{\lmod}{L-\mtr{mod}}
\newcommand{\rex}{\mtr{Rex}}
\newcommand{\resml}{\mtr{res}^A_L}
\newcommand{\cbt}{\;\underline{{\cdot}}\;}
\newcommand{\acm}{{A_{\cm}}}
\newcommand{\cfm}{\cf(\cm)}
\newcommand{\res}{{\mtr{Res}}}
\newcommand{\repk}{{\rep(K)}}
\newcommand{\resra}{{\res^{\mtr{ra}}}}
\newcommand{\rha}{{\mtr{\ra}}}
\newcommand{\epsul}{\eps_L}
\newcommand{\epsm}{\eps_M}
\newcommand{\acml}{(\al^c_M)|_L}
\newcommand{\idm}{\id_M}
\newcommand{\idl}{\id_L}
\newcommand{\deltau}{\delta_{\unu}}
\newcommand{\ptr}{\mtr{ptr}}
\newcommand{\zom}{{Z(M)}}
\newcommand{\csl}{{\cs_L}}
\newcommand{\ztmn}{{Z^{(2)}(M,N)}}
\newcommand{\bfastar}{\bfa^*}
\newcommand{\zu}{{Z(\unu)}}
\newcommand{\di}{d_i}
\newcommand{\ccro}{{\cc^{\rho}}}
\newcommand{\ccpj}{{\cc^{(j)}}}
\newcommand{\cfa}{\rc}
\newcommand{\cflo}{\cf(L_1)}
\newcommand{\cfltw}{\cf(L_2)}
\newcommand{\fjst}{F^{j}_{st}}
\newcommand{\ccpsj}{{\cc^{(j)}_s}}
\newcommand{\ccpsro}{{\cc^{\ro}_s}}
\newcommand{\ccpjs}{\ccpsj}
\newcommand{\ccjps}{\ccpsj}
\newcommand{\ccjsp}{\ccpsj}
\newcommand{\tfjst}{\widetilde{\fjst}}
\newcommand{\lj}{\mtc L_j}
\newcommand{\mtmj}{\mtc M_j}
\newcommand{\tfjpsptp}{{\widetilde{F^{j'}_{s't'}}}}
\newcommand{\tfjstp}{{\widetilde{F^{j}_{st'}}}}
\newcommand{\enxczcca}{{\enx_\czcc(A)}}
\newcommand{\ccjpt}{\cc^{(j)}_t}
\newcommand{\ccpjt}{\cc^{(j)}_t}
\newcommand{\sumjtomst}{\sum_{j\in \mtc J}\;\sum_{s,t\in \mtmj}}
\newcommand{\fjpsptp}{{F^{j'}_{s't'}}}
\newcommand{\fjstp}{{F^{j}_{st'}}}
\newcommand{\alijst}{{\al(i)^j_{st}}}
\newcommand{\csujst}{\csu^j_{st}}
\newcommand{\st}{_{st}}
\newcommand{\mtha}{\mathbb A}
\newcommand{\bcl}{\bar{\mtr{CE}}(L)}
\newcommand{\jcd}{\mtc J_{\cd}}
\newcommand{\jce}{\mtc J_{\ce}}
\newcommand{\ls}{{(L, S)}}
\newcommand{\jtwo}{j_{2}}
\newcommand{\hmda}{\hm_{D(A)}}
\newcommand{\uwha}{\uw_\lsh^A}
\newcommand{\mci}{{\mtc I}}
\newcommand{\irrda}{\irr(D(A))}
\newcommand{\lsg}{{L(g)}}
\newcommand{\lsh}{{L(h)}}
\newcommand{\dimkm}{\dimk(M)}
\newcommand{\lsm}{{L(m)}}
\newcommand{\phia}{\phi(A)}
\newcommand{\jcl}{\mtc J_L}
\newcommand{\kastar}{{K(A^*)}}
\newcommand{\cajojtw}{\ca(j_1, j_2)}
\newcommand{\mtcj}{\mtc J}
\newcommand{\cbjojtw}{\ccb(j_1, j_2)}
\newcommand{\nioitw}{N^i_{\io, \itw}}
\newcommand{\barcmio}{{\bar C}_{m(\io)}}
\newcommand{\barcmitw}{{\bar C}_{m(\itw)}}
\newcommand{\alio}{\al_{ _\io}}
\newcommand{\ali}{\al_{ _i}}
\newcommand{\alitw}{\al_{ _\itw}}
\newcommand{\barci}{\bar C_{i}}
\newcommand{\barcmi}{\bar C_{m(i)}}
\newcommand{\mio}{m(\io)}
\newcommand{\mitw}{m(\itw)}
\newcommand{\ccmio}{\cc^{m(\io)}}
\newcommand{\ccmi}{{\cc^{M(i)}}}
\newcommand{\ccmj}{\cc^{m(j)}}
\newcommand{\ccmitw}{\cc^{m(\itw)}}
\newcommand{\ccjo}{\cc^{\jo}}
\newcommand{\ccjtw}{\cc^{\jtw}}
\newcommand{\nkij}{{N^k_{ij}}}
\newcommand{\cscm}{\cs_\cm}
\newcommand{\wloge}{without loss of generality }
\newcommand{\barzmu}{ \barz_{\cm}(\unu)}
\newcommand{\reg}{\mtr{Reg}}
\newcommand{\barzrm}{ \barz_{\cm}}
\newcommand{\esuv}{F^s_{uv}}
\newcommand{\esvu}{F^s_{vu}}
\newcommand{\esuu}{F^s_{uu}}
\newcommand{\esvv}{F^s_{vv}}
\newcommand{\tildef}{\tilde F}
\newcommand{\ovl}{\overline}
\newcommand{\sccj}{{|\cc^j|}}
\newcommand{\ccl}{\cc^l}
\newcommand{\cck}{\cc^k}
\newcommand{\ccp}{{\cz_2(\cc)}}
\newcommand{\gr}{\mtr{Gr}}
\newcommand{\mtlsum}{\mathlarger{\sum}}
\newcommand{\whag}{{\widehat{(a, \gamma)}}}
\newcommand{\grcc}{\gr(\cc)}
\newcommand{\sumstom}{\sum_{s=0}^m}
\newcommand{\barcr}{\bar C_r}
\newcommand{\barcl}{\bar C_l}
\newcommand{\barcj}{\bar C_j}
\newcommand{\barck}{\bar C_k}
\newcommand{\io}{{i_1}}
\newcommand{\itw}{{i_2}}
\newcommand{\jo}{{j_1}}
\newcommand{\jtw}{{j_2}}
\newcommand{\barcjo}{\bar C_{j_1}}
\newcommand{\barcjtw}{\bar C_{j_2}} 
\newcommand{\chitw}{\chi_{\itw}}
\newcommand{\xitw}{x_{ _{\itw}}}
\newcommand{\mi}{{M(i)}}
\newcommand{\chio}{\ch_{\io}}
\newcommand{\xio}{x_{ _{\io}}}
\newcommand{\phibr}{\phi_{\bar R}}
\newcommand{\jc}{{\mtc J}^{c}}
\newcommand{\cvi}{C_{V_i}}
\newcommand{\tausg}{{\tau_{ _\sg}}}
\newcommand{\cdl}{\cd_L}
\newcommand{\cfl}{\overline{\mtr{CF}}(L)}
\newcommand{\cel}{\overline{\mtr{CE}}(L)}
\newcommand{\mtcfj}{\mtc F_j}
\newcommand{\jl}{\mtc J_L}
\newcommand{\resal}{\mtr{Res}_L}
\newcommand{\lamcdl}{\lam_{\cdl}}
\newcommand{\mmj}{\mtc M_L}
\newcommand{\zim}{\{0,\dots , m\}}
\newcommand{\zir}{\{0,\dots , r\}}
\newcommand{\cdm}{\cd_M}
\newcommand{\cdn}{\cd_N}
\newcommand{\jm}{\mtc J_M}
\newcommand{\resalfj}{\resal(F_j)}
\newcommand{\zh}{\mtc Z(H)}
\newcommand{\enxcza}{\enx_{\czcc}(A)}
\newcommand{\bcj}{\bar C_j}
\newcommand{\bpcj}{\bar C'_j}
\newcommand{\lamcdm}{\lam_{\cd_M}}
\newcommand{\lamcdn}{\lam_{\cd_N}}
\newcommand{\ena}{\enx_\czcc(A)}
\newcommand{\sumutom}{{\sum_{u=0}^m}}
\newcommand{\mti}{\mtc I}
\newcommand{\svec}{\mtr{sVec}}
\newcommand{\irrcc}{\irr(\cc)}
\newcommand{\irrcd}{{\irr(\cd)}}
\newcommand{\cocd}{\co(\cd)}
\newcommand{\cocc}{\co(\cc)}
\newcommand{\sgst}{\sg^{\mtr{st}}}
\newcommand{\sgstzv}{\sg^{\mtr{st}}_{Z(V), -}}
\newcommand{\dltv}{\delta_V}
\newcommand{\bdlt}{{\bar{\delta}}}
\newcommand{\bpivx}{\bar{\pi}_{V;\;X}}
\newcommand{\bpimx}{\bar{\pi}_{M;\;X}}
\newcommand{\bpinx}{\bar{\pi}_{N;\;X}}
\newcommand{\bpi}{\bar{\pi}}
\newcommand{\mtj}{\mtc L}
\newcommand{\mtjcd}{\mtj_{\cd}}
\newcommand{\lag}{\langle}
\newcommand{\rag}{\rangle}
\newcommand{\dvx}{\delta_{V, X}}
\newcommand{\evx}{{\ev_X}}
\newcommand{\barzw}{{\barz(W)}}
\newcommand{\bphi}{\bar{\phi}}
\newcommand{\zcat}{Z_{\mtr{cat}}}
\newcommand{\repdh}{\rep(D(H)}
\newcommand{\rrange}{R()}
\newcommand{\lrange}{L()}
\newcommand{\dimkh}{\dimk(H)}
\newcommand{\reph}{\rep_{\kk}(H)}
\newcommand{\omg}{\Omega}
\newcommand{\hdl}{\hat{\mtc D}_L}
\newcommand{\phiro}{{\phi_\ro}}
\newcommand{\tev}{\tilde{ev}}
\newcommand{\betai}{\beta^{-1}}
\newcommand{\alphai}{\alpha^{-1}}
\newcommand{\uu}{{(1,1)}}
\newcommand{\tu}{{(2,1)}}
\newcommand{\td}{{(2,2)}}
\newcommand{\ud}{{(1,2)}}
\newcommand{\dtwo}{{(2)}}
\newcommand{\bara}{{\bar{A}}}
\newcommand{\minu}{{(-1)}}
\newcommand{\zero}{{(0)}}
\newcommand{\fic}{F_i^c}
\newcommand{\chc}{\ch^c}
\newcommand{\barga}{\bar G(A)}
\newcommand{\repa}{\rep(A)}
\newcommand{\repaad}{\rep(A)_{\mtr{ad}}}
\newcommand{\repbad}{\rep(B)_{\mtr{ad}}}
\newcommand{\repb}{\rep(B)}
\newcommand{\rp}{\rep}
\newcommand{\cmhhh}{_H\cm_H^H}
\newcommand{\cmhh}{_H\cm_H}
\newcommand{\hot}{\hat{\ot}}
\newcommand{\bars}{\bar S}
\newcommand{\cmh}{_H\cm}
\newcommand{\phiu}{\phi_1}
\newcommand{\bphiu}{\bar{\phi}_1}
\newcommand{\phit}{\phi_2}
\newcommand{\bphit}{\bar{\phi}_2}
\newcommand{\phitr}{\phi_3}
\newcommand{\bphitr}{\bar{\phi}_3}
\newcommand{\qur}{q_1^{\mtr{R}}}
\newcommand{\qtwor}{q_2^{\mtr{R}}}
\newcommand{\coih}{{\mtr{co} \; H}}
\newcommand{\qhb}{\text{quasi-Hopf bimodules}\;}
\newcommand{\qtwol}{q_2^{\mtr L}}
\newcommand{\qul}{q_1^{\mtr L}}
\newcommand{\hcmhh}{\;^H_H\cm_H}
\newcommand{\rphi}{{\;_R\phi}}
\newcommand{\tblam}{{\tilde{\blam}}}
\newcommand{\fpch}{\fp(\ch)}
\newcommand{\sgi}{\sigma(i)}
\newcommand{\ccrbz}{\ccr^{\barz}}
\newcommand{\mtl}{{  \mtc L}}
\newcommand{\mtlu}{{  \mtc L}(\unu)}
\newcommand{\piuy}{\pi_{\unu, Y}}
\newcommand{\piux}{\pi_{\unu, X}}
\newcommand{\iuy}{\iota_{\unu, Y}}
\newcommand{\iux}{\iota_{\unu, X}}
\newcommand{\rhom}{\rho_M}
\newcommand{\z}{{  Z}}
\newcommand{\czx}{\cz(X)}
\newcommand{\czm}{\z(M)}
\newcommand{\czv}{\cz(V)}
\newcommand{\zm}{\cz(M)}
\newcommand{\pivx}{\pi_{V;\;X}}
\newcommand{\barza}{\barz(A)}
\newcommand{\hra}{\hookrightarrow}
\newcommand{\barj}{\mathbf J}
\newcommand{\mul}{\mu_{\barzm}^l}
\newcommand{\mur}{\mu_{\barzm}^r}
\newcommand{\chofm}{\ch_M}
\newcommand{\pium}{\pi_{\unu, M}}
\newcommand{\prc}{\perp^{\mtr{proj}}}
\newcommand{\wrt}{with respect to\;}
\newcommand{\fpas}{\frac{p_s}{a_s} }
\newcommand{\sgs}{{\sigma(s)}}
\newcommand{\ccsgs}{\cc_{\sgs}}
\newcommand{\cu}{\mtc U}
\newcommand{\czzcc}{\cz(\czcc)}
\newcommand{\cci}{\cc^i}
\newcommand{\cfzcc}{{\cf(\czcc)}}
\newcommand{\czccad}{{\czcc_{\ad}}}
\newcommand{\vsk}{\vskip 0.15cm \noindent}
\newcommand{\ccj}{\cc^j}
\newcommand{\chad}{\ch_{\ad}}
\newcommand{\sumstor}{\sum_{s=0}^r}
\newcommand{\fpcc}{\fp(\cc)}
\newcommand{\fpx}{\fp(X)}
\newcommand{\ccs}{\cc^s}
\newcommand{\cct}{\cc^t}
\newcommand{\ircc}{\irr(\cc)}
\newcommand{\ccsuu}{\ccs_{uu}}
\newcommand{\lamad}{\lam_{\ad}}
\newcommand{\regcc}{R_\cc}
\newcommand{\gch}{G_{\ch}}
\newcommand{\gmu}{G_\mu}
\newcommand{\mus}{\mu_s}
\newcommand{\gze}{G_0}
\newcommand{\chun}{\ch(1)}
\newcommand{\pkjojtw}{{\wpp_k(\jo, \jtw)}}
\newcommand{\pccj}{p_{\cc^j}}
\newcommand{\wtpccj}{\widetilde{\pccj}}
\newcommand{\bxd}{\boxed}
\newcommand{\incl}{\hookrightarrow}
\newcommand{\lha}{^A\lh}
\newcommand{\fbx}{\fbox}
\newcommand{\tsig}{\tilde\sigma}
\newcommand{\mir}{{\mtr{mir}}}
\newcommand{\wtl}{\widetilde}
\newcommand{\forg}{\mtr{Forg}}
\newcommand{\barh}{\underline H}
\newcommand{\alx}{\al_X}
\newcommand{\bxt}{\boxtimes}
\newcommand{\barhx}{\barh(X, X)}
\newcommand{\barcx}{\barc(X, X)}
\newcommand{\czcm}{\cz(\cc^*_\cm)}
\newcommand{\cstrm}{\cc^*_\cm}
\newcommand{\apk}{A//K}
\newcommand{\stc}{^{*\cop}}
\newcommand{\ka}{\mtr{K}(A)}
\newcommand{\xij}{x_{ij}}
\newcommand{\xii}{x_{ii}}
\newcommand{\barg}{{  \overline{G}}}
\newcommand{\vect}{\mtr{Vec}}
\newcommand{\jcdu}{J_{\cd}(\unu)}
\newcommand{\sent}{\mapsto}
\newcommand{\xj}{X_j}
\newcommand{\mui}{\mu_i}\newcommand{\muj}{\mu_j}
\newcommand{\kercc}{{\ker_{\cc}}}
\newcommand{\kerzcc}{{\ker_{\czcc}}}
\newcommand{\lkercc}{\lker_\cc}
\newcommand{\lkercci}{\lker_\cc(X_i)}
\newcommand{\tildefj}{\tilde{F}_j}
\newcommand{\fstar}{F_*}
\newcommand{\bargcc}{\barg(\cc)}
\newcommand{\kcc}{\mtr{K}(\cc)}
\newcommand{\sgpi}{\sg_\pi}
\newcommand{\barzx}{\barz(X)}
\newcommand{\ccrm}{\ccr_M}
\newcommand{\ccrn}{\ccr_M}
\newcommand{\barc}{{{ \ul{ \cc}}}}
\newcommand{\hmc}{\hm_\cc}
\newcommand{\rdx}{^*X}
\newcommand{\excl}{\blue{\bf !!!!!!!}}
\newcommand{\gs}{\mapsto}
\newcommand{\gos}{\mapsto}
\newcommand{\barzmn}{\barz(M\ot N)}
\newcommand{\hkw}{\hookrightarrow}
\newcommand{\qst}{q_*}
\newcommand{\qust}{q^*}\newcommand{\cec}{\mtr{  CE}(\cc)}
\newcommand{\ced}{\mtr{  CE}(\cd)}
\newcommand{\barzd}{\barz_{\cd\hookrightarrow \cc}}
\newcommand{\bll}{\blue}
\newcommand{\ovr}{\overline}
\newcommand{\pp}{\perp}
\newcommand{\dw}{\downarrow}
\newcommand{\uw}{\uparrow}
\newcommand{\dlt}{\delta}
\newcommand{\nono}{\nonumber}
\newcommand{\ch}{\chi}
\newcommand{\mtr}{\mathrm}
\newcommand{\pap}{\bowtie}
\newcommand{{\ipr}}{i'}
\newcommand{\sd}{\leq^{\oplus}}
\newcommand{\tcs}{\text}
\newcommand{\mbb}{\mathbb B}
\newcommand{\vs}{\mathbb V}
\newcommand{\sth}{suppose that\;}
\newcommand\rad{\operatorname{rad}}
\newcommand{\itm}{\item}
\newcommand{\dbd}{$$}
\newcommand{\mol}{\mtr{mod}}
\newcommand{\ro}{\rho}
\newcommand{\irr}{\mathrm{Irr}}
\newcommand{\mbc}{\mathbb C}
\newcommand{\mbs}{\mathbb S}
\newcommand{\mbz}{\mathbb Z}
\newcommand{\ct}{\mtc T}
\newcommand{\sm}{\setminus}
\newcommand{\epl}{^{+}}
\newcommand{\sbsq}{\subseteq}
\newcommand{\sbs}{\subset}
\newcommand{\cco}{\mtr{co}}
\newcommand{\cz}{\mathcal{Z}}
\newcommand{\dual}{^{*}}
\newcommand{\Gm}{\Gamma}
\ncm{\cY}{\mtc{Y}}
\newcommand{\bab}{\color{DarkOrchid}{}}
\newcommand{\eab}{\normalcolor{}}
\newcommand{\subs}{\subsection}
\newcommand{\cv}{\mtc{V}}
\newcommand{\grn}{\green}
\newcommand{\dt}{\delta}
\newcommand{\ccf}{\mathrm{ {CF}(\cc)}}
\newcommand{\cce}{\mathrm{ {CE}(\cc)}}
\newcommand{\cecc}{\mathrm{ {CE}(\cc)}}
\newcommand{\cecd}{\mathrm{ {CE}(\cd)}}
\newcommand{\kk}{\Bbbk}
\newcommand{\otL}{\ot_{L}}
\newcommand{\otl}{\ot_{L}}
\newcommand{\unpsi}{1_{\psi}}
\newcommand{\epsi}{e_{\psi}}
\newcommand{\ephi}{e_{\phi}}
\newcommand{\ech}{e_{\ch}}
\newcommand{\nleftcid}{\text{left normal  coideal subalgebra}}
\newcommand{\dimL}{\dim_{\kk}L}
\newcommand{\cl}{\mtc L}
\newcommand{\mj}{\mtc J}
\newcommand{\tl}{\tilde L}
\newcommand{\tL}{\tilde L}
\newcommand{\tpsi}{\tilde(\psi)}
\newcommand{\tmx}{\tilde{\mtc X}}
\newcommand{\zlh}{\mathrm{ZL}}
\newcommand{\ba}{\mathrm A}
\newcommand{\bv}{\mathrm V}
\newcommand{\zhopf}{\mtc{Z}_{\mtr{Hopf}}}
\newcommand{\lstar}{L^{*}}
\newcommand{\ldstar}{L^{**}}
\newcommand{\mstar}{M^{*}}
\newcommand{\mdstar}{M^{**}}
\newcommand{\lkera}{\lker_{A}}
\newcommand{\mdprime}{M''}
\newcommand{\ldprime}{L''}
\newcommand{\cm}{\mtc M}
\newcommand{\ccm}{\mathcal M}
\newcommand{\cn}{\mathcal N}
\newcommand{\ccn}{\mathcal N}
\newcommand{\rx}{\mtr{Rex}}
\newcommand{\cca}{\ca}
\newcommand{\ih}{\underline{\mtr{Hom}}}
\newcommand{\cih}{\underline{\mtr{coHom}}}
\newcommand{\hm}{\mtr{ {Hom}}}
\newcommand{\cov}{\mtr{coev}}
\newcommand{\rora}{\rho^{\mtr{ra}}}
\newcommand{\rola}{\rho^{\mtr{la}}}
\newcommand{\cx}{\mtc X}
\newcommand{\cZ}{\cz}
\newcommand{\ca}{\cA}
\newcommand{\stat}{\noindent}
\newcommand{\bfa}{{\bf A}}
\newcommand{\unu}{\mathbf{1}}
\newcommand{\barzu}{{\bar {  Z}(\unu)}}
\newcommand{\idx}{\id_X}
\newcommand{\lprime}{L'}
\newcommand{\mprime}{M'}
\newcommand{\nat}{ \mtr{{  Nat}}}
\newcommand{\ft}{\mtc F_\lam}
\newcommand{\rhau}{\rightharpoonup}
\newcommand{\lhau}{\leftharpoonup}
\newcommand{\cf}{\mathrm{ {CF}}}
\newcommand{\cfc}{\mathrm{{CF}}(\cc)}
\newcommand{\csu}{\overline{\mathfrak{  C}}}
\newcommand{\cfcc}{{\mathrm{CF}(\cc)}}
\newcommand{\catfcc}{{\mathrm{ {CF}}(\cc)}}
\newcommand{\cfcd}{{\mathrm{CF}(\cd)}}
\newcommand{\cfd}{\mathrm{CF}(\cd)}
\newcommand{\czcc}{{\cz(\cc)}}
\newcommand{\czcd}{{\cz(\cd)}}
\newcommand{\czt}{{\cz(\cz(\cc))}}
\newcommand{\rev}{\mtr{rev}}
\newcommand{\enx}{\mtr{  End}}
\newcommand{\runu}{R(\unu)}
\newcommand{\bdfn}{\bn{defn}}
\newcommand{\edfn}{\end{defn}}
\newcommand{\deltax}{\delta_X}
\newcommand{\deltav}{\delta_V}
\newcommand{\repcca}{\rep_\cc(A)}
\newcommand{\xotay}{X \ot_A Y}
\newcommand{\xoty}{X \ot Y}
\newcommand{\votw}{V \ot W}
\newcommand{\votaw}{V \ot_A W}
\newcommand{\dimax}{\dim_AX}
\newcommand{\dimccx}{\dim_\cc(X)}
\newcommand{\dimcca}{\dim_\cc(A)}
\newcommand{\dimccv}{\dim_\cc(V)}
\newcommand{\dima}{\dim_A}
\newcommand{\biga}{A}
\newcommand{\comp}{\mathbb C}
\newcommand{\tehtaa}{\theta_A}
\newcommand{\tetaa}{\theta_A}
\newcommand{\ida}{\id_A}
\newcommand{\hma}{\hm_A}
\newcommand{\hmcc}{\hm_\cc}
\newcommand{\fv}{F(V)}
\newcommand{\fw}{F(W)}
\newcommand{\ota}{\ot_A}
\newcommand{\repza}{\rep_\cc^0(A)}
\newcommand{\epsa}{\eps_A}
\newcommand{\bndefn}{\bn{defn}}
\newcommand{\edefn}{\end{defn}}
\newcommand{\bdefn}{\bn{defn}}
\newcommand{\vld}{V^{*}}
\newcommand{\vldd}{V^{**}}
\newcommand{\xld}{X^{*}}
\newcommand{\xldd}{X^{**}}
\newcommand{\yld}{Y^{*}}
\newcommand{\yldd}{Y^{**}}
\newcommand{\aldu}{A^{*}}
\newcommand{\aldd}{A^{**}}
\newcommand{\ia}{\mtr{i}_A}
\newcommand{\aota}{A\ot A}
\newcommand{\idv}{\id_V}
\newcommand{\ld}{^*}
\newcommand{\repg}{\rep(G)}
\newcommand{\thetav}{\theta_V}
\newcommand{\tta}{\theta_A}
\newcommand{\muv}{\mu_V}
\newcommand{\muw}{\mu_W}
\newcommand{\dimcc}{\dim \cc}
\newcommand{\chii}{\chi_i}
\newcommand{\chistar}{\ch_{i^*}}
\newcommand{\chj}{\ch_j}
\newcommand{\chm}{\ch_m}
\newcommand{\chn}{\ch_n}
\newcommand{\dimvi}{\dim(V_i)}
\newcommand{\mtcd}{Q}
\newcommand{\mtca}{\mtc A}
\newcommand{\lamcd}{\lam_\cd}
\newcommand{\fpdimcd}{\fp(\cd)}
\newcommand{\laml}{\lam_L}
\newcommand{\apm}{A//M}
\newcommand{\apl}{A//L}
\newcommand{\repapm}{\rep(\apm)}
\newcommand{\repapl}{\rep(\apl)}
\newcommand{\dimvj}{\dim(V_j)}
\newcommand{\dvi}{\dim(V_i)}
\newcommand{\dvj}{\dim(V_j)}
\newcommand{\sumjtom}{\sum_{j=0}^m}
\newcommand{\sumitom}{\sum_{i=0}^m}
\newcommand{\sij}{s_{ij}}
\newcommand{\sji}{s_{ji}}
\newcommand{\dxj}{d_j}
\newcommand{\dxi}{\di }
\newcommand{\dimka}{\dim_{\kk}(A)}
\newcommand{\dimk}{\dim_{\kk}}
\newcommand{\blaml}{\blam_L}
\newcommand{\sumjtor}{\sum_{j=0}^r}
\newcommand{\dimkl}{\dim_{\kk}(L)}
\newcommand{\mtcjl}{\mtc J_L}
\newcommand{\vota}{ V\ot A}
\newcommand{\vi}{V_i}
\newcommand{\vj}{V_j}
\newcommand{\dimcd}{\dim \cd}
\newcommand{\alij}{{\al_{ _{ij}}}}
\newcommand{\alji}{{\al_{ _{ji}}}}
\newcommand{\rcc}{{r_{ _{\cc}}}}
\newcommand{\rcd}{{r_{ _{\cd}}}}
\newcommand{\clsx}{[X]}
\newcommand{\clsy}{[Y]}
\newcommand{\clsz}{[Z]}
\newcommand{\rcdp}{r_{\cd'}}
\newcommand{\sumjtorp}{\sum_{j=0}^{r'}}
\newcommand{\aljm}{{\al_{ _{jm}}}}
\newcommand{\aljn}{{\al_{ _{jn}}}}
\newcommand{\sjm}{s_{jm}}
\newcommand{\smj}{s_{mj}}
\newcommand{\snj}{s_{nj}}
\newcommand{\betaij}{\beta_{ _{ij}}}
\newcommand{\betaji}{\beta_{ _{ji}}}
\newcommand{\gammaij}{\gamma_{ _{ij}}}
\newcommand{\gammaji}{\gamma_{ _{ji}}}
\newcommand{\ip}{{i'}}
\newcommand{\sumjtoprp}{\sum_{j=0}^{r'}}
\newcommand{\sumjtopr}{\sum_{j=0}^{r}}
\newcommand{\teh}{\tilde{h}}
\newcommand{\cdp}{{\cd'}}\newcommand{\xphii}{X_{\phi(i)}}
\newcommand{\inv}{^{-1}}
\newcommand{\fq}{{\mtr f_{ Q}}}
\newcommand{\tr}{\mtr{tr}}
\newcommand{\rtwone}{R_{21}R}
\newcommand{\ccad}{{\cc_{\mtr{ad}}}}
\newcommand{\ccpt}{{\cc_{\mtr{pt}}}}
\newcommand{\qtr}{quasi-triangular\;}
\newcommand{\trq}{\tr_q}
\newcommand{\repal}{\mtr{Rep}(A//L)}
\newcommand{\lkeravi}{\lker_A(V_i)}
\newcommand{\lkeravj}{\lker_A(V_j)}
\newcommand{\blml}{\blam_L} 
\newcommand{\phir}{\phi_R}
\newcommand{\kda}{{  \Phi(A)}}
\newcommand{\gpt}{{G_\pt}}
\newcommand{\mtcr}{{\mtc R}}
\newcommand{\mtcil}{\mtc{I}_L}
\newcommand{\un}{\unu}
\newcommand{\tfl}{\mtc{T}}
\newcommand{\barzm}{\barz(M)}
\newcommand{\barzn}{\barz(N)}
\newcommand{\ccr}{\mtc R^{\cc}}
\newcommand{\ulc}{\ul{\cc}}
\newcommand{\jcdpt}{{\overline{\mtc J}_{\cd_\pt}}}
\newcommand{\pimx}{\pi_{M;\;X}}
\newcommand{\pinx}{\pi_{N;\;X}}
\newcommand{\acc}{{\mathrm A_\cc}}
\newcommand{\epsu}{\eps_\unu}
\newcommand{\jpte}{{\mtcj_{\cc_\pt}^e}}
\newcommand{\ob}{\mtr{Obj}}
\newcommand{\obc}{\mtr{Obj(\cc)}}
\newcommand{\ccop}{\cc^{\mtr{op}}}
\newcommand{\mtf}{\mtc F_\lam}
\newcommand{\mtfi}{\mtc F^{-1}_\lam}
\newcommand{\wti}{{\widetilde{\mtc I}}}
\newcommand{\wtj}{{\widetilde{\mtc J}}}
\newcommand{\wkcc}{\widehat{K(\cc)}}
\newcommand{\ron}{{\frac{\ro}{x_1}}}
\newcommand{\elcd}{\ell_\cd}
\newcommand{\mcid}{\mtc I_\cd}
\newcommand{\wtlam}{\widetilde{\lam}}
\newcommand{\mcidp}{\mtc I_{\cd'}}
\newcommand{\wtildelcd}{\widetilde{\elcd}}
\newcommand{\wtildelcdp}{\widetilde{\ell_{\cd'}}}
\newcommand{\cpt}{\cc_{\mtr{pt}}}
\newcommand{\cdpte}{{\cd_\pt^e}}
\newcommand{\mtcv}{{\mtc V}}
\newcommand{\barzr}{\barz_\cd}
\newcommand{\barzv}{\barz(V)}
\newcommand{\acd}{\mathrm A_\cd}
\newcommand{\czrcd}{\cz_\cc(\cd)}
\newcommand{\sml}{\Small}
\newcommand{\bs}{\blue{\Small }}
\newcommand{\yd}{Yetter-Drinfeld}
\newcommand{\sumitor}{\sum_{i=0}^r}
\newcommand{\cdop}{\cd^{\mtr{op}}}
\newcommand{\ccrev}{\cc^{\mtr{rev}}}
\newcommand{\barz}{{\bar{\mathrm Z}}}
\newcommand{\etl}{etale\;}
\newcommand{\czca}{\cz(\ca)}
\newcommand{\clcrd}{\mtr{Cl}(\cc:\cd)}
\newcommand{\tetx}{\text}
\newcommand{\widehta}{\widehat}
\newcommand{\wdhat}{\widehat}
\newcommand{\wht}{\widehat}
\newcommand{\cofa}{{\mathbb C[\mtc B]}}
\newcommand{\wdt}{\widehat}
\newcommand{\dl}{{^\#}}
\newcommand{\comx}{\mathbb C}
\newcommand{\sgj}{{\sg(j)}}
\newcommand{\mujo}{\mu_\jo}
\newcommand{\mujtw}{\mu_\jtw}
\newcommand{\adz}{a^{\#}}
\newcommand{\bdz}{b^{\#}}
\newcommand{\spr}{S^\perp}
\newcommand{\cofs}{\comp [S]}
\newcommand{\spz}{S^{\perp_z}}
\newcommand{\omz}{\omega_z}
\newcommand{\zg}{\mathrm{Z}(S)}
\newcommand{\aling}{{\al \in g}}
\newcommand{\blkg}{\mtr{Bl}(g)}
\newcommand{\clsg}{\mtr{Cl}(g)}
\newcommand{\mtadinv}{\mtc G^{{-1}}}
\newcommand{\muk}{\mu_{k}}
\newcommand{\mta}{\mtc F}
\newcommand{\cofad}{\comp[\wdht A]}
\newcommand{\wtau}{\wdht{\tau}}
\newcommand{\mtainv}{{\mta}^{-1}}
\newcommand{\wdht}{\widehat}
\newcommand{\augm}{\mtr{aug}}
\newcommand{\mua}{\wdht {\wdht a}}
\newcommand{\aps}{A//S}
\newcommand{\krn}{\mtr{ker}}
\newcommand{\cssa}{\cc(S, A)}
\newcommand{\aug}{\mtr{aug}}
\newcommand{\rss}{{\big|_S}}
\newcommand{\gprp}{g^\perp}
\newcommand{\alins}{{s \in S}}
\newcommand{\gm}{\gamma}
\newcommand{\sz}{s^{D}}
\newcommand{\wmu}{\widehta{\mu}}
\newcommand{\wmui}{\widehta{\mu}_i}
\newcommand{\wmuj}{\widehta{\mu}_j}
\newcommand{\wch}{\widehta{\ch}}
\newcommand{\cO}{\mtc{O}}
\newcommand{\wzd}{\widehat{d}}
\newcommand{\wpm}{\widehat{P}}
\newcommand{\wps}{\widehat{p}}
\newcommand{\gal}{\mtr{Gal}}
\newcommand{\galkq}{\gal(\mathbb K/\mathbb Q)}
\newcommand{\sgh}{\sg_{ _{H}}}
\newcommand{\sggi}{{\sg(i)}}
\newcommand{\sge}{\sg_{_{\widehat R}}}
\newcommand{\unue}{{\unu_{\cecc}}}
\ncm{\rk}{\mtr{rank}}
\newcommand{\mtcf}{\mtc {F}}
\newcommand{\core}{\mtr{core}}
\newcommand{\wsgf}{\widehat{{\sg}_{ _F}}}
\newcommand{\sghstar}{{{\sg}_{ _{H^*}}}}
\newcommand{\we}{\widehta{E}}
\newcommand{\sumktom}{\sum_{k=0}^m}
\newcommand{\wf}{\widehat{F}}
\newcommand{\hsgj}{\widehat{\sg}(j)}
\newcommand{\whsgi}{\widehta{\sg}(i)}
\newcommand{\wpp}{\widehat{p}}
\newcommand{\dimccj}{\dim(\cc^j)}
\newcommand{\tauj}{{\tau(j)}}
\newcommand{\dimcctauj}{\dim(\cc^\tauj)}
\newcommand{\etas}{{\eta(s)}}
\newcommand{\mcc}{{M_\cc}}
\newcommand{\wal}{\widehta{\al}}
\newcommand{\wj}{\widehat{\mtc J}}
\newcommand{\galc}{\mtr{Gal}_{\cc}}
\newcommand{\galz}{\mtr{Gal}_{\czcc}}
\newcommand{\wjr}{\widehat{J}_{R}}
\newcommand{\dimcck}{\dim(\cc^k)}
\newcommand{\wgrcc}{\widehat{\mtr{Gr}(\cc)}}
\newcommand{\nchi}{{\frac{\ch_i}{\di }}} \newcommand{\nchj}{{\frac{\ch_j}{\dxj}}}
\newcommand{\wni}{{\widehat{n}_i}}
\newcommand{\sgte}{\widetilde{\sg_E}}
\newcommand{\mtad}{\mtc G}
\newcommand{\whj}{\widehta{h}_j}
\newcommand{\jdl}{{j\dl}}
\newcommand{\wcfcc}{\widehat{\cfcc}}
\newcommand{\mutauj}{\mu_{\tau(j)}}
\newcommand{\tauk}{\tau(k)}
\newcommand{\muzm}{{\mu_0^{-}}}
\newcommand{\sqrtog}{\sqrt{|G|}}
\newcommand{\muz}{\mu_0}
\newcommand{\njtw}{n_\jtw}
\newcommand{\njo}{n_\jo}
\newcommand{\fjo}{F_\jo}
\newcommand{\fjtw}{F_\jtw}
\newcommand{\wta}{\widehat{A}}
\newcommand{\dol}{{^{\circ}}}
\newcommand{\bdl}{{b\dl}}
\newcommand{\jdol}{{j\dol}}
\newcommand{\fj}{F_j}
\newcommand{\cwta}{\comp[\wta]}
\newcommand{\hx}{\widehta{x}}
\newcommand{\hy}{\widehta{y}}
\newcommand{\cal}{\mtc A_{\al}}
\newcommand{\xuu}{x_{uu}}
\newcommand{\wxuu}{\widehat{\xuu}}
\newcommand{\xvv}{x_{vv}}
\newcommand{\xuv}{x_{uv}}
\newcommand{\xmn}{x_{m,n}}
\newcommand{\buvmn}{B^{u,v}_{m,n}}
\newcommand{\blm}{\blam}
\newcommand{\dimccr}{\dim(\cc^r)}
\newcommand{\adl}{a\dl}
\newcommand{\sumltom}{\sum_{l=0}^m}
\newcommand{\mbq}{\mathbb Q}
\newcommand{\mbqs}{\mathbb Q(S)}
\newcommand{\mbk}{\mathbb K}
\newcommand{\mz}{\mathbb Z}
\newcommand{\wsgj}{\widehat{\sigma}(j)}
\newcommand{\wsgi}{\widehat{\sigma}(i)}
\newcommand{\wg}{\widehat{g}}
\newcommand{\wtf}{\widehat{F}}
\newcommand{\galqspq}{\mtr{Gal}(\mathbb Q(S)/\mathbb Q)}
\newcommand{\cctauj}{\cc^{\tau(j)}}
\newcommand{\cctauk}{\cc^{\tau(k)}}
\newcommand{\wtfj}{\widetilde{F_j}}
\newcommand{\wfj}{\widetilde{F_j}}
\newcommand{\wtmuj}{\widetilde{\mu_j}}
\newcommand{\wmtcfj}{{\widetilde{\mtc F}_j}}
\newcommand{\mtfr}{\mtr{F_a}}
\newcommand{\wdr}{R_\comp^*}
\newcommand{\fgph}{{F_{G/H}}}
\newcommand{\wcfj}{\wmtcfj}
\newcommand{\nxi}{{\frac{x_i}{\di }}}
\newcommand{\fpr}{{\fp(R)}}
\newcommand{\nxs}{{\frac{x_s}{d_s}}}
\newcommand{\mtfme}{\mtc F}
\newcommand{\chic}{\ch_i^{\circ}}
\newcommand{\chjc}{\ch_j^{\circ}}
\newcommand{\mtfsh}{{\mtc F_\lam}}
\newcommand{\mupq}{{\mu_{pq}}}
\newcommand{\tlam}{{\widetilde{\lam}}}
\newcommand{\chid}{{\ch_i^{\circ}}}
\newcommand{\rc}{{R_\comp}}
\newcommand{\rgo}{{\mathbb R_{\geq 0}}}
\newcommand{\sumrorc}{{{\sum\limits_{\ro \in \rc}}}}
\newcommand{\aliro}{{\ro(x_i)}}
\newcommand{\tomega}{\widetilde{\omega}}
\newcommand{\barjd}{{\bar{\mtc J_\cd}}}
\newcommand{\lbarcj}{{\frac{C_j}{{\dim(\mathcal C^j)}}}}
\newcommand{\omtcb}{{\overline{\mathcal B}}}
\newcommand{\whr}{{\widehat{R}}}
\newcommand{\nxj}{\frac{x_j}{\dxj}}
\newcommand{\nxk}{\frac{x_k}{d_k}}
\newcommand{\onkij}{{\overline{N^k_{ij}}}}
\newcommand{\sgk}{{\sigma(k)}}
\newcommand{\sgl}{{\sigma(l)}}
\newcommand{\fqi}{{\fq^{-1}}}
\newcommand{\wdb}{{\widehat{\mtc B}}}
\newcommand{\mtcb}{{\mtc B}}
\newcommand{\nif}{{h_i}}
\newcommand{\rb}{(R, \mtc B)}
\newcommand{\nxip}{{\frac{x_{i'}}{d_{i'}}}}
\newcommand{\mujp}{{\mu_{j'}}}
\newcommand{\etai}{{\eta(i)}}
\newcommand{\wsg}{{\widehat{\sg}}}
\newcommand{\wsgh}{{\wsg_{ _{H}}}}
\newcommand{\wtaui}{{\widehat{\tau}(i)}}
\newcommand{\wsghstar}{{\wsg_{H^*}}}
\newcommand{\wtauj}{{\wtau(j)}}
\newcommand{\weta}{{\widehat{\eta}}}
\newcommand{\detai}{{d_{ _{\eta(i)}}}}
\newcommand{\hbz}{{(H, \mtc B, \mu_0)}}
\newcommand{\whbz}{{\widehat{\hbz}}}
\newcommand{\tsgh}{{{\widetilde{\sgh}}}}
\newcommand{\ghb}{{G_\hb}}
\newcommand{\taujo}{{\tau(\jo)}}
\newcommand{\taujtw}{{\tau(\jtw)}}
\newcommand{\mutauk}{{\mu_{\tauk}}}
\newcommand{\hetai}{{h_{ _\etai}}}
\newcommand{\xetai}{{x_{ _\etai}}}
\newcommand{\wn}{{\widehat{n}}}
\newcommand{\wh}{{\widehat{h}}}
\newcommand{\distar}{{d_{i^*}}}
\newcommand{\dwtaui}{{d_{ _{\wtaui}}}}
\newcommand{\sumiptom}{{\sum_{\ip=0}^m}}
\newcommand{\alitaugj}{{\al_{ _{i\tau_g(j)}}}}
\newcommand{\dtaugj}{{d_{ _{\tau_g(j)}}}}
\newcommand{\mip}{{M(\ip)}}
\newcommand{\sumttom}{{\sum_{t=0}^m}}
\newcommand{\muxi}{{\mu_{ _{[X_i]}}}}
\newcommand{\muxip}{{\mu_{ _{[X_{\ip}]}}}}
\newcommand{\ncj}{{\frac{C_j}{\dim(\cc^j)}}}
\newcommand{\jp}{{j'}}
\newcommand{\minv}{{M^{-1}}}
\newcommand{\tfq}{{\widehat{\fq}}}
\newcommand{\catcecc}{{\mtr{CE}(\cc)}}
\newcommand{\wcatfcc}{{\widehat{\catfcc}}}
\newcommand{\mforall}{{\;\;\text{for all}\;\;}}
\newcommand{\what}{\widehat}
\newcommand{\wfz}{{\widehat{F}_0}}
\newcommand{\hbfr}{{(H, \mtc B, \fp)}}
\newcommand{\sgn}{{\mtr{sgn}}}
\newcommand{\wir}{{\widehat{R}}}
\newcommand{\gcc}{{G(\cc)}}
\newcommand{\jccpt}{{J_{ _{\ccpt}}}}
\newcommand{\jccad}{{J_{ _{\ccad}}}}
\newcommand{\fpccad}{{\fp(\ccad)}}
\newcommand{\fpccpt}{{\fp(\ccpt)}}
\newcommand{\kc}{{K(\cc)}}
\newcommand{\wkc}{{\widehat{\kc}}}
\newcommand{\hs}{{(L,\;\mtc S)}}
\newcommand{\rbad}{{H_{ _{ad}}}}
\newcommand{\hbad}{{\hb_{ad}}}
\newcommand{\jhbad}{{J_{\hbad}}}
\newcommand{\htt}{{(K, \mtc T)}}
\newcommand{\jhtt}{{\mtc J_{ _{\htt}}}}
\newcommand{\jhs}{{\mtc J_{ _{\hs}}}}
\newcommand{\lamhs}{{\lam_{ _{\hs}}}}
\newcommand{\lamhtt}{{\lam_{ _{\htt}}}}
\newcommand{\coo}{{co}}
\newcommand{\wdhad}{{(\wdh)_{ad}}}
\newcommand{\had}{{H_{ad}}}
\newcommand{\wdh}{\widehat{H}}
\newcommand{\whbad}{{\whb_{ _{ad}}}}
\newcommand{\kerhb}{{\ker_{ _{\hb}}}}
\newcommand{\gwdh}{{G(\wdh)}}
\newcommand{\nxl}{\frac{x_l}{d_l}}
\newcommand{\proditom}{{\prod_{i=0}^m}}
\newcommand{\prodjtom}{{\prod_{j=0}^m}}
\newcommand{\wdhn}{{\wdh^{(n)}}}
\newcommand{\hn}{{H_{(n)}}}
\newcommand{\nox}{{\fp(x)}}
\newcommand{\mujstar}{{\mu_{j^\#}}}
\newcommand{\sco}{{S^\coo}}
\newcommand{\rrad}{{I(1)}}
\newcommand{\istar}{{i^*}}
\newcommand{\mtcs}{{\mtc S}}
\newcommand{\lamccad}{{{\lam_{\ccad}}}}
\newcommand{\lamczccad}{{{\lam_{\cczc_{ad}}}}}
\newcommand{\wtm}{{\widetilde{m}}}
\newcommand{\cczerog}{{\cc^0_G}}
\newcommand{\zd}{{\mathbb Z_d}}
\newcommand{\ucc}{{U(\cc)}}
\newcommand{\uczcc}{{U(\czcc)}}
\newcommand{\czccpt}{{\czcc_{pt}}}
\newcommand{\cg}{{[G,G]}}
\newcommand{\wtp}{{\widehat{P}}}
\newcommand{\pie}{{\pi_e}}
\newcommand{\cdpt}{{\cd_{\pt}}}
\newcommand{\focc}{{\fpcc}}
\newcommand{\wtlp}{{\widehat{p}}}
\newcommand{\wlamczccad}{{\widetilde{\lam}_{\czccad}}}
\newcommand{\wlamccad}{{\widetilde{\lam}_{\ccad}}}
\newcommand{\oucc}{{|\ucc|}}
\newcommand{\ccadpt}{{{\ccad}_{pt}}}
\newcommand{\oG}{{|G|}}
\newcommand{\oN}{{|N|}}
\newcommand{\oH}{{|H|}}
\newcommand{\ccgad}{{\cc^G_{ad}}}
\newcommand{\dimCi}{{\dim(\cc^i)}}
\newcommand{\dimCj}{{\dim(\cc^j)}}
\newcommand{\dimCk}{{\dim(\cc^k)}}
\newcommand{\ckij}{{c^k_{ij}}}
\newcommand{\pkij}{{p^k_{ij}}}
\newcommand{\Vi}{{V_i}}
\newcommand{\wtd}{\widetilde}
\newcommand{\mtrce}{{\mtr{CE}}}
\newcommand{\wtmu}{\widetilde{\mu}}
\newcommand{\ceczcc}{\mtrce(\czcc)}
\newcommand{\tch}{\widetilde{\ch}}
\newcommand{\chpp}{{\ch^{''}}}
\newcommand{\lbardm}{\overline{d}_m}
\newcommand{\lbardt}{\frac{\mtr{D}_t}{\dim(\cd^t)}}
\newcommand{\xisstar}{{X_{s^*}}}
\newcommand{\chit}{\ch_t}
\newcommand{\chitstar}{\ch_{t^*}}
\newcommand{\chisstar}{\ch_{s^*}}
\newcommand{\chis}{{\ch_{ _{i_s}}}}
\newcommand{\ccpp}{{\cz_2(\cc)}}
\newcommand{\ma}{{\mtr A}}
\newcommand{\sgf}{\sg_{_F}}
\newcommand{\mtjccp}{\mtc J_{\ccp}}
\newcommand{\mtjcdp}{\mtj(\cdp)}
\newcommand{\mtjce}{\mtj(\ce)}
\newcommand{{\bm}}{{\overline{m}}}
\newcommand{\sumrtom}{\sum_{r=0}^m}
\newcommand{\nzm}{{Z(M)}}
\newcommand{\ru}{R(\unu)}
\newcommand{\ml}{{\mtr L}}
\newcommand{\fcjps}{{F(\ccjps)}}
\newcommand{\fcjpt}{{F(\ccjpt)}}
\newcommand{\tp}{{t'}}
\newcommand{\otmtha}{{\ot_\ma}}
\newcommand{\mmino}{{m_{-1}}}
\newcommand{\nmino}{{n_{-1}}}
\newcommand{\bigop}{\big(}
\newcommand{\bigcp}{\big)}
\newcommand{\dmino}{{_{-1}}}
\newcommand{\dzero}{{_{0}}}
\newcommand{\cresal}{\mtr{CFRes}^A_L}
\newcommand{\pil}{{\pi}}
\newcommand{\xic}{{x_i^{\circ}}}
\newcommand{\xjc}{{x_j^{\circ}}}
\newcommand{\dimccpt}{{\dim(\ccpt)}}
\newcommand{\piz}{\pi_e}
\newcommand{\kcd}{{K(\cd)}}
\newcommand{\jccp}{{\mtc J_{ _\ccp}}}
\newcommand{\xit}{X_t}
\newcommand{\xitstar}{{X_{t^*}}}
\newcommand{\xis}{X_s}
\newcommand{\rj}{{\mtc R_j}}
\newcommand{\wtpres}{{{\wtp{\big|_{\ccp}}}}}
\newcommand{\fpccp}{{\fp(\ccp)}}
\newcommand{\fprj}{{\fp(\rj)}}
\newcommand{\jcc}{{\mtcj_{\cc}}}
\newcommand{\jcdp}{{\mtcj_{\cdp}}}
\newcommand{\fpccj}{{\fp(\ccj)}}
\newcommand{\dimccp}{{\dim(\ccp)}}
\newcommand{\dimrj}{{\dim(\rj)}}
\newcommand{\stabj}{{\mtr{stab}_j}}
\newcommand{\cccd}{{{\small(\cc\big/ \cd)_r}}}
\newcommand{\fpcdpt}{\fp(\cdpt)}
\newcommand{\vpt}{{\mtc P^{\pt}}}
\newcommand{\qpt}{{\mtc Q^{\pt}}}
\newcommand{\vf}{{\mtc P^{f}}}
\newcommand{\qf}{{\mtc Q^{f}}}
\newcommand{\vp}{{\mtc P}}
\newcommand{\fpczccpt}{{\fp(\czccpt)}}
\newcommand{\vpccpt}{{\vp(\ccpt)}}
\newcommand{\alp}{{\al_p}}
\newcommand{\cccsd}{{(A/B)_r}}
\newcommand{\mh}{\mtc H}
\newcommand{\almcd}{\lamcd} 
\newcommand{\xim}{x_m}
\newcommand{\is}{{\mtc I_S}}
\newcommand{\mcci}{{M_{\cc}(i)}}
\newcommand{\mcdi}{{M_{\cd}(i)}}
\newcommand{\wpkij}{{\wp_k(i,j)}}
\newcommand{\dimcci}{{\dim(\cc^i)}}
\newcommand{\nchk}{{\frac{\ch_k}{d_k}}}
\newcommand{\jinjccp}{{j\in \jccp}}
\newcommand{\jinjccpt}{{j\in \jccpt}}
\newcommand{\jinjcd}{{j\in \jcd}}
\newcommand{\ellcd}{{\ell_\cd}}
\newcommand{\classnu}{{\mtr{Class}(\nu)}}
\newcommand{\mcd}{{M_\cd}}
\newcommand{\wtc}{{\widehat{\mtr{C}}}}
\newcommand{\fqcd}{{\fq^\cd}}
\newcommand{\fqcc}{{\fq^\cc}}
\newcommand{\rjcapcd}{{\mtc {R}(\cd)_j}}
\newcommand{\mzcd}{{\cz_2(\cd}}
\newcommand{\brj}{{\overline{\rj}}}
\newcommand{\bdr}{\bn{dr}}
\newcommand{\edr}{\end{dr}}
\newcommand{\bnq}{\bn{quest}}
\newcommand{\enq}{\end{quest}}
\newcommand{\bnnd}{\bn{nd}}
\newcommand{\ennd}{\end{nd}}
\newcommand{\fpccg}{{\fp(\cc_{g})}}
\newcommand{\azp}{{b_0}}
\newcommand{\azq}{c_0}
\newcommand{\azpq}{{a^0_{pq}}}
\newcommand{\agp}{b_g}
\newcommand{\agq}{c_g}
\newcommand{\agpq}{{a^g_{pq}}}
\newcommand{\irrcdp}{{\irr(\cdp)}}
\newcommand{\aly}{\al_Y}
\newcommand{\alz}{\al_Z}
\newcommand{\cdppt}{{\cdp_\pt}}
\newcommand{\ccg}{{\cc_g}}
\newcommand{\psixi}{{\psi_{[X_i]}}}
\newcommand{\psixip}{{\psi_{[X_\ip]}}}
\newcommand{\lamcdp}{{\lam_\cdp}}
\newcommand{\stabc}{{\mtr{stab}}}
\newcommand{\ccptcapccad}{{\ccpt\cap\ccpt}}
\newcommand{\cdad}{{\cd_{\ad}}}
\newcommand{\fpcdp}{{\fp(\cdp)}}
\newcommand{\sone}{{s_1}}
\newcommand{\alprime}{{\al^{'}}}
\newcommand{\hacd}{{\mh(\cd)}}
\newcommand{\kcrcd}{{\mh(\cd)}}
\newcommand{\ccrcp}{{K(\cc:\ccp)}}
\newcommand{\lamccp}{\lam_{ _\ccp}}
\newcommand{\crcdad}{{\ccrcd_{\ad}}}
\newcommand{\ucrcdad}{{U_{\ccrcd_{\ad}}}}
\newcommand{\ukcrcd}{{U_{\kcrcd}}}
\newcommand{\ofp}{{\overline{\fp}}}
\newcommand{\rcdxi}{{R_{ _\cd}(X_i)}}
\newcommand{\rcdxj}{{R_\cd(X_j)}}
\newcommand{\rcdxk}{{R_\cd(X_k)}}
\newcommand{\wtlamcd}{{\wtlam_\cd}}
\newcommand{\kcad}{{K(\ccad)}}
\newcommand{\rnu}{{r_\nu}}
\newcommand{\brnu}{{R_\nu}}
\newcommand{\kcrcdad}{{\kcrcd_{\ad}}}
\newcommand{\pkcrcd}{{P_{ _{\kcrcd}}}}
\newcommand{\xnu}{{X_\nu}}
\newcommand{\pcd}{{P_{ _\kcrcd}}}
\newcommand{\pcdsq}{{P_{ _\kcrcd}^2}}
\newcommand{\fco}{{co}}
\newcommand{\cde}{{\cd^E}}
\newcommand{\stabcdx}{{\stabc_{\cd}(X)}}
\newcommand{\stcd}{{\stabc_\cd}}
\newcommand{\expp}{{\exp_p}}
\newcommand{\fpxt}{{\fp(X_t)}}
\newcommand{\fpxs}{{\fp(X_s)}}
\newcommand{\ckad}{{\ck_{ad}}}
\newcommand{\oucaad}{{U(\ccad)}}
\newcommand{\caad}{{\ca_{ad}}}
\newcommand{\fpcaad}{{\fp(\caad)}}
\newcommand{\capt}{{\ca_{pt}}}
\newcommand{\fpcapt}{{\fp(\capt)}}
\newcommand{\caadpt}{{(\ca_{ad})_{\pt}}}
\newcommand{\fpcaadpt}{{\fp(\caadpt)}}
\newcommand{\caadad}{{{\ca^{(2)}_{ad}}}}
\newcommand{\caadadpt}{{\ca^{(2)}_{ad}}_{\pt}}
\newcommand{\xp}{{X_p}}
\newcommand{\ucd}{{U(\cd)}}
\newcommand{\rg}{{\mtr{rank}(\cc_g)}}
\newcommand{\barcp}{{\overline{P}_{ _\cc}}}
\newcommand{\barcpg}{{\overline{P}_{ _\ccg}}}
\newcommand{\barcpone}{{\overline{P}_{ _{\cc_1}}}}
\newcommand{\rctwo}{{r_{\cc^{(2)}}}}
\newcommand{\wkcaveecb}{{\widehat{\mk(\caveecb)}}}
\newcommand{\mk}{{\mtr{K}}}
\newcommand{\wkcd}{{\widehat{\mk(\cd)}}}
\newcommand{\wkca}{{\widehat{\mk(\ca)}}}
\newcommand{\wkcb}{{\widehat{\mk(\ccb)}}}
\newcommand{\kca}{{\mk(\ca)}}
\newcommand{\kcb}{{\mk(\ccb)}}
\newcommand{\caveecb}{{\ca\vee\ccb}}
\newcommand{\DB}{\text{dual-Burnside}}
\newcommand{\db}{\text{dual-Burnside}}
\newcommand{\cbad}{{\ccb_{ad}}}
\newcommand{\lamcdad}{{\lam_\cdad}}
\newcommand{\cfca}{{\cf(\ca)}}
\newcommand{\cfcb}{{\cf(\cb)}}
\newcommand{\kccp}{{K(\ccp)}}
\newcommand{\rankcc}{{\rank(\cc)}}
\newcommand{\scrg}{{\mathscr G}}
\newcommand{\scra}{{\mathscr A}}
\newcommand{\scrb}{{\mathscr B}}
\newcommand{\scrc}{{\mathscr C}}
\newcommand{\scrd}{{\mathscr D}}
\newcommand{\wnu}{{\widetilde{\nu}}}
\newcommand{\nx}{{\frac{[X]}{\fp(X)}}}
\newcommand{\ny}{{\frac{[Y]}{\fp(Y)}}}
\newcommand{\gwkc}{{G(\wkc)}}
\newcommand{\gwkd}{{G(\wkd)}}
\newcommand{\wkd}{{\widehat{K(\cd)}}}
\newcommand{\chx}{{\ch_X}}
\newcommand{\wirrc}{{\widehat{\irrcc}}}
\newcommand{\genx}{{\langle X\rangle}}
\newcommand{\zccx}{{Z_\cc(X)}}
\newcommand{\hlrd}{{\mtr{Cl}({H}:{L})}}
\newcommand{\pu}{\green{P-U}\;}
\newcommand{\pij}{{\pi(j)}}
\newcommand{\mux}{{\mu_X}}
\newcommand{\ji}{{j_i}}
\newcommand{\muji}{{\mu_{j_i}}}
\newcommand{\zcdx}{{Z_\cd(X)}}
\newcommand{\roc}{{\ro_{ _C}}}
\newcommand{\kzc}{{K_0(\cc)}}
\newcommand{\wirrcd}{{\widehat{\irr(\cd)}}}
\newcommand{\oucd}{{|U(\cd)|}}
\newcommand{\dimx}{{\dim(X)}}
\newcommand{\ndim}{{n_{dim}}}
\begin{document} 

\title[Isaacs type results for fusion categories]
{On some integral properties of 
dimensions in Isaacs fusion categories}

\begin{abstract}
For a fusion category, we prove some new integral properties concerning the dimension of a simple object that generates a Isaacs fusion subcategory. If any such fusion subcategory is Isaacs this implies that $\cc$ is a Frobenius type fusion category. A stronger divisibility result is proven for any modular fusion category. This divisibility result implies the  converse direction of a Ito-Michler type result for modular fusion categories,  see \cite[Proposition 2.3]{b-conj-im}. 
\end{abstract}

\thanks{This work was supported by a grant of the Ministry of Research, Innovation and Digitization, CNCS -UEFISCDI, project number PN-IV-P1-PCE-2023-0264, within PNCDI IV}

\author{Sebastian Burciu}
\address{Inst.\ of Math.\ ``Simion Stoilow" of the Romanian Academy P.O. Box 1-764, RO-014700, Bucharest, Romania}
\email{sebastian.burciu@imar.ro}
\maketitle
\setcounter{tocdepth}{1}
\tableofcontents
\section{Introduction}
Fusion categories represent a significant advancement in both mathematics and theoretical physics, providing a categorical framework for exploring symmetry and structure across diverse contexts, including quantum field theory, topology, and representation theory. By generalizing classical concepts from group theory, fusion categories address various settings  where traditional group symmetries fail to describe the intricacies of complex phenomena.

A central problem in the theory of fusion categories is the generalized Kaplansky's 6th conjecture, which states that in any spherical fusion category, the dimensions of all simple objects divide the dimension of the category itself. 

This conjecture plays a crucial role in the study of fusion categories by offering insights into the structural properties of these mathematical objects.  The conjecture is significant because it establishes a deeper connection between fusion categories and representation theory. Moreover, verifying Kaplansky's sixth conjecture could lead to advancements in understanding the interplay between quantum symmetries and topological features, further enriching the theoretical landscape of both mathematics and physics.

In this paper, we establish new divisibility results for the dimensions of simple objects in fusion categories, generalizing classical results from finite group representation theory. For instance, it is well-known that for a finite group $G$, the dimension of any of its complex irreducible representations divides the index of the center of the group, see \cite[Theorem 3.12]{is}. We show that this property extends to a broad class of  fusion categories:
\bt\label{ribbon-uc}
Let $\cc$ be a psudo-unitary ribbon category, and $X$ be a simple object of $\cc$.  Then,
$$
\fp(X)\mid \frac{\fp(\cc)}{\oucc}.
$$
\et
This result is a special case of a more general divisibility criterion, presented in Theorem \ref{isaacs-gf-uc-zccx}), involving the dimension of  the {\it center} of a simple object (defined in Equation \eqref{z-cent}).

For modular fusion categories, we prove a stronger constraint:
\bt\label{modular-uc}
Let $\cc$ be a  modular fusion category, and let $X$ be a simple object of $\cc$ such that $\cc=\genx$. Then,
$$
\dim(X)^2\mid \frac{\dimcc}{\oucc}.
$$
\et
These theorems arise naturally in the study of $s$-Isaacs fusion categories, a concept introduced in \cite{eno-nec} to investigate integrality properties related to the Liu-Palcoux-Wu invariants and the generalized Kaplansky's 6th conjecture. This concept builds upon the notion of Isaacs fusion rings as discussed in \cite{lpr1, lpr2}. 

Recall \cite{eno-nec}, that a fusion category $\cc$ is called {\it $s$-Isaacs} if for any character $\muj: \kc\to \comp$ and any simple object $X\in \cc$, the quantity 
$$
\lambda_s(\muj,X):=(\dim \cc)^s\dimccj^{1-s}\frac{\muj([X])}{\dim(X)}
$$
is an algebraic integer. For the definition of the conjugacy class $\cc^j$ corresponding to the character $\mu_j$, see Subsection \ref{bij-conj}.

The two theorems presented above are particular cases of the following more general results concerning $s$-Isaacs fusion categories.

\bt \label{s-geq-1/2-gf-uc}
Let $s\geq 1/2$ be a rational number and $\cc$ be a  spherical fusion category such that $\wkc$ is real non-negative (RN). Suppose that $X$ is a simple object such that the fusion subcategory $\genx$ is $s$-Isaacs. Then:
$$
\frac{(\dimcc)^{2s}\fpcc}{\dim(X)^2\oucc}\in \mathbb A.
$$
\et
The above theorem follows from a more general version, see Theorem \ref{s-geq-1/2-gf-zccx}, involving the aforementioned notion of the center of a simple object. The definition of the dual space $\wkc$ and  the {\it real non-negative (RN)} condition are described in Section \ref{prelim}.
\bt \label{s-geq-0-gf-uc}
Let $s\geq 0$ be a rational number and $\cc$ be a  spherical fusion category such that $\wkc$ is real non-negative (RN). Suppose that $X$ is a simple object such that the fusion subcategory $\genx$ is $s$-Isaacs. Then:
$$
\frac{(\dimcc)^{s}\fpcc}{\dim(X)\oucc}\in \mathbb A.
$$
\et 
The above theorem also generalize further (see Theorem \ref{s-geq-0-gf-zcx},) by incorporating \emph{the center} of a simple object. Notably, since every ribbon category is Isaacs ($s$ = 0) \cite[Proposition 5.2]{eno-nec}, and  $\wkc$ is real non-negative (RN) in this case, these results apply broadly to ribbon fusion categories.

\subsection{On the converse part of Ito-Michler}
A celebrated result in finite group theory, the Ito-Michler theorem, see \cite{michler}, states that for a finite group $G$ and a prime  $p$ dividing the order of $G$ the following are equivalent:
\bne [label=(\alph*)]
\item
No irreducible character of  $G$ has degree divisible by $p$.
\item
A Sylow $p$-subgroup of $G$ is both abelian and normal. 
\ene

In \cite{b-conj-im} we explored an analogue of $(a)\implies (b)$ for modular categories. 
More precisely, we have proved the following proposition:
\bp\label{mtc-deg-p-not div}\cite[Proposition 2.3]{b-conj-im}\\
Let $\cc$ be a weakly integral modular fusion category, and let $p$ be a prime divisor of $\fpcc$ such that $\fpcc=p^\al N$ with $\al\geq 1$ and $(p, N)=1$. If $p$ does not divides the Frobenius-Perron dimension  of  any of its simple objects,  then $p^{\al}$ divides $|U(\cc)|$, the order of the universal grading group of $\cc$.  
\ep
In this paper, we use Theorem \ref{modular-uc} to demonstrate the other direction, $(b)\implies (a)$.  Notably, this direction can be extended to the slightly broader context of braided categories:
\bt\label{im-converse}
Suppose that $\cc$ is a weakly integral braided fusion category. If $\fp(\cc)=p^\al N$ with $(p,N)=1$ and $p^{\al}\mid \oucc$, then $\fp(X)^2$ is not divisible by $p$ for any simple object $X$ of $\cc$.
\et
Note that since $\cc$  is weakly integral then $\fpx^2\in \mathbb Z$ for all any $X\in \irrcc$, see \cite[Theorem 3. 10]{ng}.

For modular categories, this yields a full analogue of the Ito-Michler theorem:
\bc[Complete Ito-Michler analogue]\label{complete-im}
Let $\cc$ be a weakly integral modular fusion category, and let $p$ be a prime divisor of $\fpcc$ such that $\fpcc=p^\al N$ with $\al\geq 1$ and $(p, N)=1$. Then $p$ does not divides the Frobenius-Perron dimension  of  any of its simple objects  if and only if $p^{\al}$ divides $|U(\cc)|$, the order of the universal grading group of $\cc$.  
\ec

The organization of the paper is as follows:

Section \ref{prelim} reviews basic facts concerning pivotal fusion categories from \cite{scalg} and proves a key divisibility result,  see Theorem \ref{uc-div},  concerning the size of the center of a simple object.
Section \ref{bantay} generalize a character class identity from \cite[Proposition 5.1]{b-conj-im}.  In Section \ref{s-geq 1/2} we prove Theorem \ref{s-geq-1/2-gf-uc}. In section \ref{s-geq-0} we prove Theorem  \ref{s-geq-0-gf-uc}. In the last Section we prove Theorem \ref{im-converse} and study a consequence of \cite[Conjecture 7.2]{eno-nec}.
\section{Preliminaries}\label{prelim}
Throughout this paper, we work over the base field $\mathbb{C}$.  For $\alpha, \beta \in \mathbb{C}$ with $\beta \neq 0$, we write $\beta \mid \alpha$ if $\alpha/\beta \in \mathbb{A}$, the ring of algebraic integers.

For a detailed exposition on fusion categories, we recommend the monograph \cite{EGNO15}.

Let $\cc$ be a fusion category.  We will refer to the maximal pointed fusion subcategory of $\cc$ as $\ccpt$, and let $\irr(\cc)$ represent a set of representatives for the isomorphism classes of simple objects within $\cc$. The cardinality of $\irr(\cc)$ is called the \emph{rank of $\cc$},  denoted as $\rank(\cc)$. We let $\kc:=\kzc\otimes_{\mathbb Z} \comp$, where $\kzc$ is the Grothendieck ring of $\cc$.

\subsection{Pivotal and Spherical Structures}
A \emph{pivotal structure} on a rigid monoidal category $\cc$ is a monoidal isomorphism $j:\id_\cc\to ()^{**}$ between monoidal functors. A fusion category with such a structure is called \emph{pivotal}.  As in \cite{EGNO15}, a pivotal fusion category is  \emph{spherical} if and only if $\dim(V ) = \dim(V^*)$ for all simple objects $V$, where $\dim$ denotes the \emph{quantum dimension} of $\cc$. A fusion category $\cc$ is \emph{pseudo-unitary} if $\fp(\cc)=\dimcc$. In this context, according to \cite[Proposition 8.23]{eno-annals}, the category $\cc$ possesses a unique spherical structure, designated as the \emph{canonical spherical structure}, which ensures that the categorical dimensions of all simple objects are positive. Under this spherical structure, the categorical dimension of any object $X$ in $\cc$ coincides with its Frobenius–Perron dimension, i.e., $\fp(X) = \dim(X)$ holds for every object $X$. A fusion category is referred to as \emph{weakly integral} if its Frobenius-Perron dimension is an integer. Moreover, it is termed \emph{integral} if the Frobenius-Perron dimension of each object is an integer. By \cite[Proposition 8.24]{eno-annals}, weakly integral fusion categories are pseudo-unitary.

\subsection{Drinfeld center and adjunction}
For any fusion category $\cc$, its Drinfeld center $\mathcal{Z}(\cc)$  is a braided fusion category $\czcc$ equipped with a forgetful functor $F:\czcc\ra \cc$, as described in \cite{EGNO15}.

The forgetful functor $F$ has a right adjoint functor denoted as $R:\cc \ra \czcc$, allowing us to define a Hopf comonad $Z :=FR:\cc \ra \cc$ as an end.
  
Indeed, following \cite[Section 2.6]{scalg}, it is established that  
\beq
Z(V)\simeq \int_{X\in \cc}X\ot V\ot X^*
\eeq

It is known that the object $A:=Z(\unu)$ (the image of the unit object) constitutes a (central) commutative algebra $\cz(\cc)$, referred to as \emph{the adjoint algebra} of $\cc$.

The vector space $\cecc:= \hm_{\cc}(\unu, A) $ is called the {\it the set of central elements}. The vector space  $\cfcc:=\hm_\cc(A,\unu)$ is known as {\it the space of class functions} of $\cc$.

For a  fusion category over  $\comp$, according to \cite[Example 4.4]{scalg}, there is a $\comp$-algebra isomorphism $\cfcc\simeq \kzc\ot_{\Z}\comp$. Therefore, $\cfcc$ is commutative if and only if $\kzc$ is commutative. Both  $\cfcc$ and $\cecc$ possess a semisimple $\comp$-algebra structure, as described in \cite{scalg}.  Furthermore,  $\cecc$ is a commutative algebra, implying that it decomposes as a direct product of copies of the field $\comp$.

Given that $R:\cc \ra \czcc$ is a right adjoint to the forgetful functor $F:\czcc \ra \cc$,   \cite[Theorem 3.8]{scalg} establishes that this adjunction induces an isomorphism of $\comp$-algebras:
\beq\label{adjisom}
\cfcc \xra{\cong} \mtr{End}_{\czcc}(R(\unu)),\;\; \ch\mapsto Z (\ch)\circ \delta_\unu.
\eeq
where $\delta_\unu$ represents the comultiplication of the Hopf comonad $Z$ associated at the unit $\unu$ of the fusion category $\cc$. 
\bn{exmp}
If $H$ is a semisimple Hopf algebra, then, for $\cc=\rep(H)$, it follows from \cite[Subsect. 3.7]{scalg} that $\cfcc\simeq C(H)$, the \emph{character algebra} of $H$. Recall that $C(H)$ consists of all linear trace maps $f:H\ra \comp$.
\end{exmp}
\subsection{A bijection on the set of conjugacy class sums} \label{bij-conj} 

Assume that $\cc$ is a pivotal fusion category with a commutative Grothendieck ring $K(\cc)$. Since $\czcc$ is also a fusion category, we can express $R(\unu)=\bigoplus_{i=0}^m\mathcal C_i$ as a direct sum of simple objects in $\czcc$. It is noteworthy that any two distinct simple objects $\cc^i$ are not isomorphic to one another, since from the isomorphism of Equation \eqref{adjisom} the endomorphism ring of $R(\unu)$ is a commutative $\comp$-algebra. We refer to,  $\mathcal C^{0},\dots, \mathcal C^{m}$ as {\it the conjugacy classes of $\cc$. } Since the unit object $\unu_{\czcc }$ is clearly a sub-object of $R(\unu)$, we can assume $\mathcal C^{0} = \unu_{\czcc }$.

In \cite[Section 5]{scalg}, Shimizu also  established a canonical pairing between $\cfcc$ and $\cecc$ with values in $\comp$,  defined by the relation:
\beq\label{pairing-def}
f\circ a=\langle f, a\rangle\id_\unu
\eeq
for all $f \in \cfcc$ and $a\in \cecc$.

Let $\tilde F_0, \tilde F_1, \dots, \tilde F_m\in \enx_{\czcc}(R(\unu))$ denote the canonical projections onto each of the conjugacy classes. The corresponding primitive idempotents of $\cfcc$ under the canonical adjunction isomorphism $\cfcc \simeq \enx_{\czcc}(R(\unu))$ from Equation \eqref{adjisom} will be denoted as $F_0, F_1, \dots , F_m$. This establishes a canonical  bijection between the (central) primitive idempotents of $\cfcc$ and the conjugacy classes of $\cc$. Let $\mtcj:=\{0,1,\dots ,m\}$ represent the set of indices corresponding to the central primitive idempotents of  $\cfcc$. For each $F_j$ with $j\in \mtc J$, the aforementioned bijection assigns a conjugacy class $\ccj$. Since $\kc$ is commutative, we have $\dim_\comp\kc=m+1$ and $\{F_j\}$ for a $\comp$-linear basis of $\kc$.

\subsection{On the definition of conjugacy class sums}
In a pivotal fusion category, Shimizu defines the elements $C_j := \mtc{F}_\lam^{-1}(F_j) \in \cecc$ as the \emph{conjugacy class sums} corresponding to the conjugacy class $\mathcal{C}^j$. Here, $\lambda \in \cf(\cc)$ represents a cointegral of $\cc$ such that $\langle \lambda, u \rangle = 1$ (refer to \cite[Section 5]{scalg}). Moreover, the morphism  $u:\unu\ra A$ represents the unit of the central algebra $A$, and $\mtc F_{\lambda}$ is the \emph{Fourier transform} of $\cc$ associated with $\lambda$.  This is  a  linear map defined by:
\beq
\mtc F_{\lambda}:\cecc\ra \cfcc\;\;\text{given by}\;\;a \mapsto \lambda \lh \mtc S(a)
\eeq
where $\mtc{S}$ denotes the antipode.

\subsection{The integral of a fusion category}
Let $\cc$ be a fusion category, and let $A=Z(\unu)$ be its adjoint algebra as previously defined. In accordance with  \cite[Definition 5.7]{scalg}, an \emph{ integral} in $\cc$ is represented as a morphism $\Lambda: \unu \ra A$ in $\cc$ satisfying
\beqn
 m \circ (\id_{A}\otimes \Lambda)=\eps_{\unu}\ot \Lam,
\eeqn
where $m:A\otimes A\ra \unu$ denotes the multiplication of the adjoint algebra $A$. Here, $\eps_\unu$ represents the counit $\eps_V:Z(V)\ra V $ specialised at $V=\unu$, the unit object of the fusion category $\cc$.

It is established that the integral of a  fusion category is unique up to a scalar.  Let $\Lam\in \cecc$ denote the idempotent integral of $\cc$ and define $\widehat{\Lam}:=\dim(\cc)\Lam$. 
 
Let $\{V_0, \dots, V_m\}$ be a complete set of representatives of the isomorphism classes of simple objects in $\cc$. Shimizu provided, in \cite[Section 6]{scalg},  a construction of the (central) primitive elements $E_i\in \cecc$,  corresponding to the irreducible characters $\ch_i:=\ch_{V_i}$ associated to each simple object $V_i$ of $\cc$. This correspondence is characterized by the condition: $$\lag\ch_i, E_j\rag=\delta_{i, j},$$ for all $i,j$.  

Assuming that $\kc$ is commutative, note that  \cite[Lemma 4.1]{ccc-march}, along with Equation (4.7) from the same paper, provides an alternative pair of dual bases for the canonical pairing. More explicitly, it follows that $\{F_j\}$ and $\{\frac{C_j}{\dim(\ccj)}\}$ form another pair of dual bases expressed as:
\beq\label{dbase}
\langle F_j, \frac{C_{j'}}{\dim(\cc^{\jp})} \rangle=\delta_{j,j'}
\eeq
for all $j, j'\in \mtc J$.

Without loss of generality, we can set $V_0=\unu$, the unit of the fusion category $\cc$. Consequently, we have $\ch_0=\epsu$, and the associated primitive idempotent is $E_0=\Lam$, as seen in \cite[Lemma 6.1]{scalg}.

By \cite[Lemma 2.1]{b-conj-im} we have that
\beq\label{integrform}
\widehat{\Lam}=\sum_{j\in \mtcj} C_j.
\eeq 
for any  a pivotal fusion category $\cc$ with a commutative Grothendieck ring.
\subsection{Orthogonality relations}
Let $\cc$ be a fusion category with a commutative  Grothendieck ring $\kzc$.  Define $\wirrc$ as the set of all linear characters $\mu: \kc \to \comp$. Given that $\kc$ is commutative, $\wirrc$ forms a basis for the linear dual space $\kc^*$, denoted by $\wkc$.  Thus, we have $\wkc:=\comp[\wirrc]$.

In $\kc$, the orthogonality relations are as follows. The first orthogonality relation from \cite[Equation (2.5)]{b-pa} expressed as:
\beq\label{first:gen:orth}
\sumitom \mu_j([X_i])\mu_k([X_{i^*}])=\delta_{j,k} n_j,
\eeq
where $n_j$ represents the \emph{formal codegree} of $\muj$, as as described in \cite[Subsection 2.1]{b-pa} or  \cite{O3}. For $j=k$, this relation simplifies to:
\beq\label{nj}
\sumitom |\mu_j([X_i])|^2=n_j.
\eeq
The second orthogonality relation in $\kc$,  from \cite[Equation (2.6)]{b-pa}, (see also \cite[Corollary 6.11]{scalg}), can be written as:
\beq\label{second:orth-2}
\sumjtom\frac{\dimccj}{\dimcc}\muj([X])\muj([Y^*])=\delta_{[X], [Y]},
\eeq
for any simple objects $X$ and $Y$ of $\cc$. 
\subsection{$\wkc$ multiplication} 
Let $\cc$ be a pivotal fusion category with Grothendieck group  $\kzc$. Recall $\kc=\kzc\otimes_{\mathbb Z} \comp$. We define a multiplication on $\wkc$ as follows.  For two linear characters $\mui, \muj \in \wirrc$, their multiplication is given by:
\beq\label{mwa}
[\mu_i\star \mu_j](\nx):=\mui(\nx)\muj(\nx), \;\text{for all}\;X\in \irr(\cc).
\eeq
This multiplication can be extended linearly to the entire space $\wkc$, yielding an algebra structure on $\wkc$. 

Additionally, there exists an involution on $\wkcc$ defined on the basis $\wirrc$ by $\muj \mapsto \mu_{j^\#}$, where $\mu_{j^\#}$ is defined by $\mu_{j^\#}([X]) := \muj([X^*])$ for all $X \in \irrcc$. It is also well-known that $\muj([X^*])=\overline{\muj([X])}$ for all $X\in \irrcc$.

Without loss of generality, we may suppose that $\mu_0:=\fp$, i.e $\mu_0([X])=\fpx$, for all $X\in \irrcc$. Equation \eqref{mwa} implies that $\mu_0$ is the unit of $\wkc$. Following the terminology in \cite[Section 2]{b-pa}, $\wkc$ becomes the dual hypergroup of the abelian normalizable hypergroup $\kc$. Since $\wirrc$ is a linear basis for $\wkc$, there exist non-zero scalars ${\widehat p}_k(i,j)\in \comp$ such that:
\beq\label{hatp:eq}
\mui \star \muj=\sumktom{\wdht p}_k(i,j)\muk.
\eeq 
We refer to ${\wdht p}_k(i,j)$ as \emph{dual-structure fusion coefficients.}. A formula for the dual fusion coefficients ${\wdht p}_k(i,j)$  is provided in \cite[Proposition 2.1]{b-blms}. In the fusion category settings, this  can be expressed as:
\beq\label{hat:pk-2}
\widehat{p}_k(iweakly integral,j)=\frac{1}{n_k}\bigg(\sum_{X\in \irrcc}\frac{1}{\fpx}\mu_i([X])\mu_j( [X])\mu_k([X^*])\bigg),
\eeq
where $n_k$, as above, it is the formal codegree of $\mu_k$.  According to  \cite[Equation (4.7)]{ccc-march}, for any spherical fusion category, we have
\beq\label{nj-sph}
n_j=\frac{\dim(\cc)}{\dim(\cc^j)}.
\eeq  
See also \cite[Theorem 2.13]{O3} for a related statement.

For a weakly integral fusion category $\cc$, by \cite[Proposition 8.27]{eno-annals} it follows that $\dim(\cc^j)\in \mathbb Z$.
\subsection{Kernels and centers of objects}
Let $\cc$ be a fusion category with a commutative Grothendieck ring. For any object $X\in \irr(\cc)$, as in \cite{bmonat}, we define its {\it kernel} as follows, 
\beq
\ker_{\cc}(X):=\{\psi\in \wirrc\;|\;\psi([X])=\fpx\}.
\eeq
The \emph{center of an object} $X \in \cc$ is defined as:
\beq\label{z-cent}
Z_\cc(X)=\{\psi\in \wirrc\;|\; |\psi([X])|=\fpx\}.
\eeq
Following \cite[Definition 5.2]{b-pa}, we can also define the dual notion of kernels as:
\beq\label{ker-psi}
\ker_{\wkcc}(\psi):=\{X\in \irr(\cc)\:|\;\psi(X)=\fpx\}.
\eeq
According to \cite[Lemma 5.3]{b-pa}, the set $\ker_{\wkcc}(\psi)$ is closed under tensor products and taking duals, thereby generating a fusion subcategory of $\cc$. Additionally, we define the dual notion of the center as:\beq\label{z-psi}
Z_{ _\wkc}(\psi)=\{X\in \irrcc\;|\; |\psi(X)|=\fpx\}.
\eeq
By \cite[Lemma 5.3]{b-pa},  the set $Z_{\wkcc}(\psi)$ is also closed under tensor products and taking duals, thus generating a fusion subcategory of $\cc$.
\subsection{Group-like elements of $\wkc$} 
Following \cite[Definition 1.3]{b-pa}, we say that an algebra morphism $
\muj\in \wkc$ is a \emph{group-like element} if $\muj\star \mu_{j^\#}=\widehat{p}_0(j, j^\#)\mu_0$. Evaluating both sides of the above equation at the unit $[\unu]$ of $\kc$, we find that $\widehat{p}_0(j, j^\#)=1$. Therefore,  $\muj\in \wkc$ is a group-like element of $\wkc$ if and only if 
\beq\label{def-g}
\muj\star \mu_{j^\#}=\mu_0.
\eeq
\br\label{gwkc-uc}
Since $\kc$ is an real non-negative (RN) hypergroup, Theorem \cite[Theorem 7.8]{b-pa} implies that for any fusion category $\cc$ we have $U(\cc)\simeq \gwkc$. 
\er
We now prove the following:
\bl
Let $\cc$ be a fusion category with a commutative Grothendieck ring. With the above notations, we have:
$$
\gwkc=\bigcap_{X \in \irrcc} \zccx
$$
\el
\bpf
For every $\mu \in \wirrc$, as mentioned in \cite[Lemma 2.21]{b-pa}, it follows that $\mu \in \gwkc$ if and only if the condition below is satisfied for all $X \in \irrcc$:
\beq\label{int-center}
|\mu([X])|=\fp(X).
\eeq 
This establishes that $\gwkc = \bigcap_{X \in \irrcc} \zccx$.
\epf
\subsection{Faithful objects and kernels of objects}
Following the definition from \cite{bmonat}, we say that an object $X \in \irrcc$ is  \textit{faithful} if the fusion subcategory subcategory generated by $X$, denoted $\genx$,  coincides with the entire category $\cc$.

An analogue of  Brauer's theorem for fusion categories with commutative Grothendieck rings was established in  \cite[Theorem 3.9]{bmonat},  and further generalized in \cite[Theorem 5.6]{b-pa}. According to this theorem, an object $X$ of $\cc$ is faithful if and only if $\kercc(X)$ is trivial.

For any simple object of $X$, we define a character of $\wkc$:
$$
\omega_X:\wkc \ra \comp,\; \omega_X(\muj)=\muj(\frac{[X]}{\fp(X)}).
$$
This character is referred to as the \emph{central character} associated with $X$ in \cite{b-pa}. From the multiplication defined on $\wkc$ in Equation \eqref{mwa}, it follows that $\omega_X$ is a an algebra homomorphism, thereby restricting to a group homomorphism:
$
\omega_X\big|:\gwkc\ra\comp^*.
$
\bl\label{faith-center}
Let $\cc$ be a fusion category with a commutative Grothendieck group. If $X$ is a faithful simple object of $\cc$, then $\omega_X$ is a faithful character of $\gwkc$.
\el
\bpf
This statement follows directly from the fusion category analogue of Brauer's theorem,  as presented  in \cite[Theorem 5.6]{b-pa}. Indeed, if $\omega_X(\mu)=1$, then  it follows that  $\mu\in \kercc(X)$, which implies that $\mu=\mu_0$.
\epf
\subsection{Orders and Perpendicular Categories}
For any linear character of $\kc$, $\muj:\kc\ra \comp$, we define its {\it order} by $\wh_j:=\widehat{p}_0(j, j^\#)^{-1}$.

\bn{defn}
For any subset $\cs\subseteq \wirrc$ we define: 
\beq\label{n-def}
n(\mtc \cs):=\sum_{\muj\in \mtc S}\whj.
\eeq
\end{defn}
According to \cite[Equation (2.10)]{b-blms} (see also \cite[Equation (2.10)]{b-pa}), we have:
\beq\label{whj:eq}
\wh_j=\frac{\fpcc}{n_j}.
\eeq 
For \textit{spherical} fusion categories, as noted in \cite[Remark 9.1]{b-pa}, Equation \eqref{nj-sph} implies:
\beq\label{whj}
\wh_j=\frac{\fpcc}{n_j}=\frac{\fpcc}{\dimcc}\dim(\ccj),
\eeq
see also \cite[Theorem 2.13]{O3}. 

Equation \eqref{whj} implies that for pseudo-unitary fusion categories, we have
\beq\label{whj-pu}
\whj=\dim(\cc^j),
\eeq
see also \cite[Equation (9.2)]{b-pa}.
\bl\label{gkcd}
Let $\cc=\genx$ be a fusion category generated by a simple object $X$ of $\cc$. Then $\zccx=G(\wkc)$.
\el
\bpf
If $\mu \in \zccx$, then by \cite[Lemma 2.35]{b-pa} it follows that $|\mu(Y)|=\fp(Y)$ for any simple object $Y$ that is a constituent of some tensor power $X^{\ot m}$ of $X$. Since $\cc=\genx$, it follows that $|\mu(Y)|=\fp(Y)$ for any simple object $Y$ of $\cc$. Thus, \cite[Theorem 6.5]{b-pa} implies that $\mu\in \gwkc$. The converse follows from the same theorem, corroborated  also with \cite[Lemma 2.31]{b-pa}.
\epf
Let $\cc$ be a fusion category with a commutative Grothendieck ring $\kc$.
\bn{defn}
We say that $\wkc$ has real non-negative coefficients, or simply that \emph{$\wkc$ is (RN)}, if all  the dual-structure fusion coefficients $\widehat{p}_k(i,j)$, defined by Equation \eqref{hatp:eq}, are real non-negative.
\end{defn}
For any subset $\cs\subseteq \wirrc$ we denote
$$\cs^\perp:=\bigcap_{\psi\in \cs}\ker_\cc(\psi). $$

By \cite[Lemma 5.3]{b-pa}, (see also \cite[Lemma 3.6]{bmonat}), it follows that $\cs^\perp\subseteq \irrcc$ is a set of simple objects closed under tensor products and taking duals.  Therefore, this generates a fusion subcategory $\cd$ of $\cc$ with $\irr(\cd)=\cs^\perp$. By abuse of notation,  we also denote this fusion subcategory by $\cs^\perp$. 

Moreover, from \cite[Proposition 2.11]{hdk}, if $\wkc$ is real non-negative (RN), then $\kc$ can be normalized to a dualizable probability group, and we have:
\beq\label{n-perp}
\fp(\cs^\perp)=\frac{\fpcc}{n(\cs)},
\eeq
where $n(\cs)$ is the order of $\cs$ as defined in Equation \eqref{n-def}.

Dually, for any subset $\ca\subseteq \irrcc$ we define 
$$
\ca^\perp:=\bigcap_{X\in \irr(\ca)}\ker_\cc(X).
$$
By the same \cite[Proposition 2.11]{hdk}, if $\wkc$ is real non-negative (RN), we have 
\beq\label{d-perp}(\cs^\perp)^\perp=\cs.
\eeq 
\subsection{Relation with the adjoint subcategory}
Next lemma is an analogue of item e)  of \cite[Lemma 2.27]{is}.
\bp\label{isaacs}
Let $\cc$ be a fusion category such that. (e.g, pseudo-unitary braided, or unitary). Then for any object $X$, we have that:  
$$
Z_\cc(X)\subseteq(\genx_{ad})^\perp.
$$ 
Moreover, suppose that $X$ is a simple object of $\cc$, then 
\beq\label{isaacs-z}
\zccx=(\genx_{ad})^\perp.
\eeq
\ep
\bpf
By definition, we have:
$$
\zccx=\{\mu\in \wirrc\;|\;\big|\mu(X)\big|=\fp(X)\}.
$$
Let $\mu \in \zccx$.  The fusion subcategory  $\genx_{ad}$ is  generated by $Y\otimes Y^*$ for $Y$ a simple object of $\genx$. If $Y\in \genx$, then $Y$ is a constituent of some tensor power $X^{\ot m}$ of $X$.  Thus, by \cite[Lemma 2.35]{b-pa}, for any $\mu\in \zccx$,  we have $\mu([Y])=\omega\fp(Y)$ for some root of unity $\omega$. This implies that $\mu([Y\ot Y^*])=\fp(Y)^2$, showing that $\mu\in (\genx_{ad})^\perp$. This proves  the inclusion $Z_\cc(X)\subseteq (\genx_{ad})^\perp$. 

Now, suppose $X$ is  simple. If $\mu\in (\genx_{ad})^\perp$, then $\mu([X\otimes X^*])=\fp(X)^2$. Since $X$ is a simple object we have $\mu([X^*])=\overline{\mu([X])}$, and this implies $\big|\mu(X)\big|=\fpx$, i.e., $\mu \in \zccx$. Thus the converse inclusion also holds when $X$ is a simple object.
\epf
\bt\label{uc-div}
Let $\cc$ be a fusion category such that $\wkc$ is real non-negative (RN) and $X$ be a simple object of $\cc$. Then 
\beq\label{zccx-oucc-div}
\frac{n(\zccx)}{\oucc}\in \mathbb A.
\eeq
\et
\bpf
As noted previously, we have $\gwkc\subseteq \zccx$ and this implies
\beq\label{incl-zccx}
\zccx^\perp\subseteq \gwkc^\perp.
\eeq
Since $\kc$ is real non-negative (RN),  as explained in \cite[Subsection 6.1]{b-pa}, it follows that $\zccx^\perp$ is closed under tensor products and taking duals. Consequently,  $\zccx^\perp$ is a fusion subcategory of $\cc$. 

Furthermore, by \cite[Corollary 7.15]{b-pa}, we have $\gwkc^\perp=\ccad$. Therefore, the inclusion from Equation \eqref{incl-zccx} is an inclusion of fusion subcategories, and we can conclude that: $\fp(\zccx^\perp)\mid \fp(\ccad)$. From this, using Equation \eqref{n-perp}, we deduce that  $\frac{\fpcc}{n(\zccx)}\mid \frac{\fpcc}{\oucc}$, completing the proof.
\epf
\br\label{br-unitary}
From the proof of \cite[Theorem 1.2]{b-blms}, specifically, using the formula from \cite[Equation (4.5)]{b-blms}, it follows that  a \emph{pseudo-unitary} braided fusion category $\cc$ has non-negative coefficients $p_k(i,j)$, i.e., $\wkc$ is real non-negative (RN).

Furthermore, for any unitary fusion category $\mathcal{C}$, $\wkc$ is also real non-negative (RN), by \cite{lpw,eno-nec}. By Theorem \ref{uc-div} it follows that \beq\label{br-zccx-oucc-div}
\frac{n(\zccx)}{\oucc}\in \mathbb A.
\eeq
for any simple object $X$ of any pseudo-unitary braided (or just unitary) fusion category $\cc$. 
\er
\section{On the size of character classes}\label{bantay}
Throughout this section, let $\cc$ be a fusion category with a commutative Grothendieck ring, $K(\cc)$. Then as above, $\cfcc$ is a commutative semisimple algebra isomorphic to $\kc\cong K_0(\cc)\ot_{\Z}\comp$. As above, Let $\mtc J^{\cc}:= \{0, \dots, m\}$ represent the set of indices corresponding to the conjugacy classes of $\cc$. It is clear that the rank of $\cc$ is $\rank(\cc) = |\mtc J^\cc| = m + 1$.

Suppose $\cd$ is a fusion subcategory of $\cc$. Note that the Grothendieck ring $K(\cd)$ is also commutative, and, as in the case above, $K(\cd)\simeq \cf(\cd)$ as algebras.  We let $\mtc J^{\cd}:= \{0, \dots, m'\}$ represent the set of indices corresponding to the conjugacy classes of $\cd$. As above, it is clear that the rank of $\cd$ is $\rank(\cd) = |\mtc J^\cd| = m' + 1$.

Let $\{F_j\}_{0 \leq j \leq m}$ be  the (central) primitive idempotents of $\cfcc$, and let $\{\cc^j\}_{0\leq j \leq m}$ be their associated conjugacy classes. For each (central) primitive  idempotent $F_j\in \cfcc$, let $\muj:\cfcc\ra \comp$ denote its associated linear character. Similarly, let $\{G_t\}_{0 \leq t \leq m'}$ be the (central) primitive idempotents of  $\cf(\cd)$, with  $\{\cd^t\}_{0\leq t \leq m'}$ being  their associated conjugacy classes.  Additionally, let $\nu_t: \cfcd \to \comp$ denote the characters associated with these idempotents.

We say that two characters $\muj, \mu_{j'}:\cfcc\ra \comp$ are \emph{$\cd$-equivalent} ($\mu_j \simeq_\cd \mu_{j'}$) if their restrictions to $\cf(\cd)$ coincide.

The set of equivalence classes is denoted by $\clcrd$. There is a natural bijection between this set and the characters of $\cfcd$.
 
Moreover, an element $C\in  \clcrd$ consists of all the indices $j$ such that the restriction $\mu_j\dw_{\cfcd}$ corresponds to the same character in $\cfcd$, denoted as $\ro_C:=\nu_s$ for some $s\leq m'$. Thus, we use $\ro_C$ to represent the common restriction to  $\cfcd$ of all characters $\mu_j$ belonging to some class $C\in \clcrd$. 

As noted above, the set of character classes of $\cc$ relative to $\cd$ corresponds bijectively to the characters of $\cfcd$. Consequently,  the central primitive idempotents of $\cfcd$ are also indexed by the classes $C\in \clcrd$. We denote the central primitive idempotent associated with a class $C \in \clcrd$ by $G_C \in \cf(\cd)$. Thus, the elements $G_C$ for $C \in \clcrd$ form a linear basis for $\cf(\cd)$, consisting of the complete set of (central) primitive idempotents of the commutative $\comp$-algebra $\cf(\cd)$. For the remainder of this section, we identify the set of indices $\mtc J^\cd$ of the central primitive idempotents of $\cf(\cd)$ with the elements of $\clcrd$. For any primitive idempotent $G_C$ of $\cfcd$, as described in Subsection \ref{bij-conj}, let  $\mtr{D}_C\in \cecd$ denote  the corresponding conjugacy class sum of $\cd$ associated with $G_C$. 

In \cite[Appendix]{ccc-march} we have constructed  two $\comp$-algebra maps, denoted by
$
\pif:\cfcd \ra \cfcc \;\text{and},\;\pie:\cecc \ra \cecd.
$ 
Furthermore, by \cite[Lemma 7.7]{ccc-march} it is known that $\pif$ is a monomorphism and $\pie$ is an epimorphism.
\bn{exmp}
If $\pi:H\ra K$ is a surjective Hopf algebra homomorphism, then clearly we have $\rep(K)\subseteq \rep(H)$. In this case, $\pi_e=\pi\big|_{Z(H)}:Z(H)\ra Z(K)$ is the restriction of $\pi$ to the center $Z(H)$ of $H$. Furthermore, $\pi_f$ corresponds to the  canonical inclusion  of the character algebra $C(K)$ of $K$ into that $H$, i.e, $C(K)\hookrightarrow C(H)$.
\end{exmp}

By \cite[Theorem 3.3]{b-conj-im} for any character class $C\in \clcrd$ and any index  $j\in C$ we have:
\beq\label{piz-form}
\piz(\lbarcj)=\frac{\mtr{D}_C}{\dim(\cd^C)}
\eeq
where $D_C$ is the conjugacy class sum associated to $\ro_C$.
Following \cite[Definition 2.5]{b-pa}, for any subset $\mtc Z\in \wirrc$ we define its \emph{dimensional order} as 
\beq\label{n-dim-def}
\ndim(\mtc Z)=\sum_{\muj\in \mtc Z}\dim(\ccj).
\eeq
Next Proposition is a generalization of \cite[Prop. 5.1]{b-conj-im}.
\bp \label{indep-prop}
Let $\cc$ be a pivotal fusion category with a commutative Grothendieck ring  $K(\cc)$ and let $\cd$ be a fusion subcategory of $\cc$.  Then for any $C \in \clcrd$ we have:
\beq\label{ext-f}
\ndim(C) = \frac{\dimcc}{\dimcd}.
\eeq
\ep
\bpf
According to Equation \eqref{piz-form}, for any 
$C \in \clcrd$, there exists a conjugacy class sum $D_C\in \cecd$ such that we have: 
$
\pi_e(C_j)= \frac{\dim(\ccj)}{\dim(\cd^C)} D_C, \;\text{for all indices}\;j\in C.
$
On the other hand, \cite[Lemma 3.4]{b-conj-im}, gives that:
\beq\label{wtlam}
\pi_e( \widehat{\Lam}_\cc) =\frac{\dimcc}{\dimcd} \widehat{\Lam}_\cd.
\eeq
Recall that $\widehat{\Lam}_\cd=\sum_{C\in \clcrd} D_C.$ by  \cite[Proposition 4.1]{ccc-march}.  The same proposition applied for $\cd=\cc$ gives that 
$\widehat{\Lam}_\cc=\sum_{j\in \mtc J} C_j.$
Counting the multiplicity of $D_C$ in the left hand side of Equation \eqref{wtlam} we find:
$$
\ndim(C)=\sum_{l \in C} \dim(\cc^l) = \frac{\dimcc}{\dimcd}.
$$
\epf
\br\label{n-vs-ndim}
Recall that for any subset $\mtc Z\subseteq \wirrc$ we have defined its order as $n(\mtc Z)=\sum_{\muj\in C}\whj$.  Using Equation \eqref{whj}, for a spherical fusion category $\cc$, dimensional and formal orders relate via: 
$$
\ndim(\mtc Z)=\frac{\dimcc}{\fpcc}n(\mtc Z),
$$
which can be written as:
\beq\label{ndim}
\frac{\fpcc}{n(\mtc Z)}=\frac{\dimcc}{\ndim(\mtc Z)}.
\eeq
Thus, for a spherical fusion category,  Equation \eqref{ext-f} implies that
\beq\label{ext-f-3}
n(C) = \frac{\fpcc}{\dimcd},
\eeq
for any  character class $C \in \clcrd$
\er
\subsection{Main Size Theorem}
\bt\label{bantay-pass}
Let $\cc$ be a pivotal fusion category with a  commutative Grothendieck ring $\kc$,  and  $\cd=\genx$ with $X\in \irrcc$. Then, with the above notations, we have:
\beq\label{eq-bantay}
\frac{\dimcc}{\ndim(Z_\cc(X))}=\frac{\dimcd}{|U(\cd)|}.
\eeq
\et
\bpf
It is straightforward to see that $\zccx$ is a union of character classes $C\in \clcrd$. Specifically, if $\mu_1\simeq \mu_2$ and $\mu_1\in \zccx$, then  since $X$ is a simple object of  $\cd$, we have $\mu_2([X])=\mu_1([X])$. Consequently, $|\mu_2([X])|=\fpx$, which implies that $\mu_2\in \zccx$ as well.

Since $\cd$ is a fusion subcategory of $\cc$, there is an algebra inclusion $K(\cd)\subseteq K(\cc)$, which gives rise to a well-defined restriction map:
\beqn
\res:\wkc\ra \wkcd, \;\mu\sent \mu\big|_\wkcd.
\eeqn

Moreover, from the definition of the centers, we have $\mu\in \zccx$ if and only if $\res(\mu)\in \zcdx$. Suppose that $\zccx=C_1\sqcup C_2\dots\sqcup C_r$ with $C_i\in \clcrd$. Then, we can express the total order as: $\ndim(\zccx)=\sum_{i=1}^r\ndim(C_i)$. By Equation \eqref{ext-f}, we know that  $\ndim(C_i)=\frac{\dimcc}{\dimcd}$, for all $1\leq i\leq r$.  Therefore, we have: $\ndim(\zccx)=r \frac{\dimcc}{\dimcd}$.

On the other hand, note that $r$ is the cardinality of the set $Z_\cd(X)$. By Lemma \ref{gkcd} we have $r=|G(\wkcd)|$, and then $r=|U(\cd)|$ by Remark \ref{gwkc-uc}. This completes the proof. 
\epf
\br
Note that $n(\gwkc)=|\gwkc|$ since $\wh_j=1$ for all $\muj\in \gwkc$. For a spherical fusion category,  the above Theorem and Equation \eqref{ndim} gives  that:
\beq\label{eq-bantay-2}
\frac{\fpcc}{ n(Z_\cc(X))}=\frac{\dimcc}{\ndim(Z_\cc(X))}=\frac{\dimcd}{|U(\cd)|}.
\eeq
\er
\section{Action of dual group-like elements}
Let $\cc$ be a fusion category with a commutative Grothendieck group  $\kzc$.  By \cite[Corollary 2.39]{b-pa}, since $\kc$ is an abelian real non-negative (RN) hypergroup, it follows that $\gwkc$ is a group. Moreover, for any $\mu\in \gwkc$, left multiplication by $\mu$ permutes the dual basis $\wirrc$. Specifically, if $\mu_i\in \wirrc$ and $\mu_z\in \gwkc$ then $\mu_z\star \mu_i=\mu_j$ for some $\mu_j\in \wirrc$. 
\subsection{Invariance properties under the group action}
We begin by establishing several key invariance properties of characters under the $\gwkc$-action.
\bl\label{abs-value-ct}
Under the $\gwkc$-action on $\wirrc$ by left multiplication, for any simple object $Y\in \irrcc$, the absolute value
$\big |\mu_j([Y]) \big |$ is constant on each orbit of this action.
\el
\bpf
Suppose  $\mujo=\mujtw\star \mu_z$ for some $z\in Z$. Then we compute:
\begin{eqnarray*}
\mu_\jtw([Y])&=&[\mu_\jo\star \mu_z]([Y])=\fp(Y)[\mu_\jo\star \mu_z](\frac{[Y]}{\fp(Y)})
\\&=&\fp(Y)\mu_\jo(\frac{[Y]}{\fp(Y)}) \mu_z(\frac{[Y]}{\fp(Y)})
\\&=&\mu_\jo([Y]) \mu_z(\frac{[Y]}{\fp(Y)}).
\end{eqnarray*}
Taking absolute values and applying Equation \eqref{int-center} yields $|\mu_{j_1}([Y])| = |\mu_{j_2}([Y])|$.
\epf
\begin{lemma}\label{nj-value-ct}
Under the $\gwkc$-action, the formal codegree $n_j$ is constant on each orbit.
\end{lemma}
\bpf
Suppose that $\mu_j=\mu_z \star \mu_k$ for some $\mu_z \in \gwkc$. Since $|\mu_z(Y)|=\fp(Y)$ for all  $Y\in \irrcc$ we have:
$$
n_j=\sum_{Y\in \irrcc}|\muj(Y)|^2= \sum_{Y\in \irrcc}|\muk(Y)|^2=n_k.
$$
Thus, characters in the same orbit share the same formal codegree.
\epf

\begin{corollary}\label{conj-dim-ct}
For a spherical fusion category $\cc$, the dimension $\dim(\cc^j)$ is constant on each orbit of the $\gwkc$-action.
\end{corollary}

\begin{proof}
This follows immediately from Lemma \ref{nj-value-ct} and Equation \eqref{nj-sph}, which relates formal codegrees to dimensions in the spherical case.
\end{proof}
\subsection{Single Generator Subcategory Case $\cc = \genx$}
We now specialize to the case where $\cc$ is generated by a single simple object $X$.
\begin{defn}
A character $\mu \in \wirrc$ is \emph{non-vanishing} if $\mu([X]) \neq 0$. An orbit $\mathcal{O} \subset \wirrc$ is called a \textit{non-vanishing orbit} if it consists of non-vanishing characters.
\end{defn}

\bl\label{free-act}
Let $\cc = \genx$ with $X \in \irrcc$. Then the $\gwkc$-action on $\wirrc$ is free when restricted to any non-vanishing orbit. 
\el
\bpf
Suppose $\mu_z \star \mu = \mu$ for some non-vanishing $\mu \in \wirrc$ and $\mu_z \in \gwkc$. Evaluating at $\nx$ gives that $\mu_z(\nx)\mu(\nx)=\mu(\nx)$.  Since $\mu(\nx) \neq 0$, we conclude $\mu_z(\nx) = 1$. As $X$ is faithful, Lemma \ref{faith-center} implies $\mu_z = \mu_1$, the identity of $\gwkc$.
\epf
\section{Isaacs fusion categories with $s\geq 1/2$}\label{s-geq 1/2}
In this section we prove some results for $s$-Isaacs fusion categories for a rational number, $s\geq 1/2$. In particular, we prove Theorem \ref{s-geq-1/2-gf-uc}.
\subsection{Definition and Basic Properties}

Let $s\geq 0$ be a rational number.   
Following \cite{eno-nec}, we have the following definition:
\begin{defn}\label{sdef} 
We say that $\cc$ is $s$-{\bf Isaacs} if for any linear character $\muj: \kc\to \comp$ and any simple object $X\in \cc$,
$$
\lambda_s(\muj,X):=(\dim \cc)^s\dimccj^{1-s}\frac{\muj([X])}{\dim(X)} \in \mathbb{A}
$$
where $\mathbb{A}$ denotes the algebraic integers.
\end{defn}   
\br
Since $\dimccj$ divides $\dim \cc$, if $\cc$ is $s$-Isaacs then it is $t$-Isaacs for any $t>s$.  As in \cite{eno-nec}, the $0$-Isaacs property will simply be called the {\it Isaacs property}. This was introduced in \cite{lpr1, lpr2} and by \cite[Proposition 5.2]{eno-nec}, any ribbon fusion category is Isaacs.
\er
Recall that $\cc$ is said to be Frobenius type if dimensions of its simple objects divide $\dimcc$.
The Kaplansky 6th conjecture for fusion categories states that any spherical fusion category
is Frobenius type. 

Following \cite{eno-nec}, we say that a category $\cc$ is {\it $s$-Frobenius type} for a
rational number $s \geq 1$ if the dimensions of its simple objects divide $(\dimcc)^s$. 
Clearly, $s$-Frobenius type implies $t$-Frobenius type for any $t > s$.
\subsection{Single Generator Case}
\bl\label{s-geq-1/2-sgl} 
Suppose that $\cc=\genx$ is a spherical fusion category with a commutative Grothendieck ring generated by a simple object $X$. If $\cc$ is s-Isaacs ($s\geq 1/2$), then 
$$
\frac{(\dimcc)^{2s+1}}{|\ucc|\dim(X)^2}\in \mathbb A.
$$
\el
\bpf
Let $\{\mathcal{A}_i\}_{i=0}^r$ be the non-vanishing orbits of $\widehat{K(\mathcal{C})}$ under $G(\widehat{K(\mathcal{C})})$-action. Since the action of $\gwkc$ is free on these orbits, it follows that $|\ca_i|=|\gwkc|$ for all $0\leq i\leq r$.

From Equation \eqref{second:orth-2} we derive:
\beq\label{second:orth-cc}
\dimcc=\sum_{j=0}^m\dim(\ccj)\mu_j([X])\mu_{j^\#}([X])
\eeq
Multiplying both sides by $\frac{(\dimcc)^{2s}}{\dim(X)^2}$ we obtain
{\Small
\beq\label{second:orth-cc-s1}
\frac{(\dimcc)^{2s+1}}{\dimx^2}=\sum_{j=0}^m\dim(\cc)^s\dim(\ccj)^{1-s}\frac{\mu_j([X])}{\dimx}\dim(\ccj)^{1-s}\dim(\cc)^s\frac{\mu_{j^\#}([X])}{\dimx}\dim(\ccj)^{2s-1}
\eeq
}
and this can be written as
\beq\label{second:orth-cc-s2}
\frac{(\dimcc)^{2s+1}}{\dimx^2}=\sum_{j=0}^m\lam_s(\muj, X)\lam_s(\mu_{j^\#}, X)\dim(\ccj)^{2s-1}
\eeq
Lemma \ref{abs-value-ct} together with Corollary \ref{conj-dim-ct} imply that each summand
$$
\lambda_s(\muj,X)\lambda_s(\mu_{j^\#},X)\dim(\ccj)^{2s-1}
$$
is constant within any non-vanishing orbit $\ca_i$. Thus the above Equation can be rewritten as: 
\begin{eqnarray*}
\frac{(\dimcc)^{2s+1}}{\dimx^2}
&=&
\sum_{i=0}^r\sum_{j\in \ca_i}\lambda_s(\muj,X)\lambda_s(\mu_{j^\#},X)\dim(\ccj)^{2s-1}
\\&=&
\sum_{i=0}^r |\ca_i|
\lambda_s(\mu_\ji,X)\lambda_s(\mu_{\ji^\#},X)\dim(\cc^{\ji})^{2s-1}
\\&=& 
|\gwkc|\sum_{i=0}^r 
\lambda_s(\mu_\ji,X)\lambda_s(\mu_{\ji^\#},X)\dim(\cc^{\ji})^{2s-1},
\end{eqnarray*}
where $j_i\in \ca_i$ is arbitrarily chosen. Since  $\gwkc\simeq \ucc$ by Remark \ref{gwkc-uc} we obtain that
\beqn
\frac{(\dimcc)^{2s+1}}{\oucc\dimx^2}=\sum_{i=0}^r 
\lambda_s(\mu_\ji,X)\lambda_s(\mu_{\ji^\#},X)\dim(\cc^{\ji})^{2s-1}.
\eeqn
Since $s\geq 1/2$ and $\dim(\ccj)\in \mathbb A$  it follows that $\dim(\ccj)^{2s-1}\in \mathbb A$. Since $\cc$ is s-Isaacs we also have $\lambda_s(\mu_\ji,X),\;\lambda_s(\mu_{\ji^\#},X) \in \mathbb A$ and the proof is completed. 
\epf
\br
If $\cc$ is spherical and $s$-Isaacs  for $s\geq 1/2$ then $\cc$ is $1+s/2$ Frobenius, see \cite[Proposition 5.4]{eno-nec}. Note that previous lemma  strengthen this result in the case of an Isaacs fusion category generated by a simple object. More precisely, under these hypotheses we derive that:
$$
\frac{(\dimcc)^{s+1/2}}{\dim(X)\sqrt{\oucc}}\in \mathbb A.
$$

\er
\subsection{General Case for $s\geq 1/2$.}
\bt \label{s-geq-1/2-gf-zccx}
Let $s\geq 1/2$ be a rational number and $\cc$ be a  spherical fusion category. Suppose that $X$ is a simple object of $\cc$ such that the fusion subcategory $\genx$ generated by $X$ is $s$-Isaacs. Then 
$$
\frac{(\dimcc)^{2s+1}}{\dim(X)^2\ndim(\zccx)}\in \mathbb A.
$$
\et
\bpf
Let $\cd:=\genx$. Applying Lemma \ref{s-geq-1/2-sgl}  for $\cd$ we have that:
$$
\dim(X)^2\mid \frac{(\dimcd)^{2s+1}}{\oucd }.
$$
By Equation \eqref{eq-bantay} we have
\beq\label{eq-bantay-pr}
\frac{(\dimcd)^{2s+1}}{\oucd }=\frac{\dimcd}{|U(\cd)|}(\dimcd)^{2s}=\frac{\dimcc}{\ndim(Z_\cc(X))}(\dimcd)^{2s}
\eeq
and therefore $$
\dim(X)^2\mid \frac{(\dimcc)^{2s+1}}{\ndim(Z_\cc(X))},
$$
since $\dim(\mathcal{C})/\dim(\mathcal{D}) \in \mathbb{A}$, see \cite{EGNO15}.
\epf
\subsubsection{Proof of Theorem \ref{s-geq-1/2-gf-uc}}
\bpf
Theorem \ref{s-geq-1/2-gf-zccx}  establishes
that $\dim(X)$ divides $\frac{(\dimcc)^{2s+1}}{\ndim(\zccx)}$. Using Equation \eqref{ndim} we obtain that $\dim(X)$ divides
\beq\label{eq-bantay-22}
(\dimcc)^{2s}\frac{\dimcc}{\ndim(\zccx)}=(\dimcc)^{2s}\frac{\fpcc}{n(\zccx)}
\eeq
On the other hand, the RN property of $\widehat{K(\mathcal{C})}$ gives that $|U(\mathcal{C})| \mid n(Z_\mathcal{C}(X))$, by Theorem \ref{uc-div}. Then the required divisibility follows.
\epf
\section{Isaacs fusion categories with $s\geq 0$}\label{s-geq-0}
In this section, we establish divisibility results for $s$-Isaacs fusion categories when $s \geq 0$, culminating in the proof of Theorem \ref{s-geq-0-gf-uc}.
\subsection{Single Generator Case}
\bl\label{gen-x}
Let $\cc=\genx$ be a spherical fusion category with a commutative Grothendieck ring, generated by a simple object $X$. Then:
\beq\label{eq-act}
\frac{\dimcc}{|\gwkc|}=\sum_{i=0}^r\dim(\cc^{j_i})\mu_{j_i}(X)\overline{\mu_{j_i}([X])},
\eeq
where $\mu_{j_i}\in \ca_i$ is arbitrarily chosen. 
\el
\bpf
Let $\ca_0, \ca_1, \dots, \ca_r\subseteq \wirrc$ be the non-vanishing orbits of the aforementioned group action. Since  $\gwkc$ actss freely on these orbits, we have that $|\ca_i|=|\gwkc|$ for all $0\leq i\leq r$.

Starting from Equation \eqref{second:orth-2}, which we rewrite as:
\beq\label{second:orth-cc}
\dimcc=\sum_{j=0}^m\dim(\ccj)\big|\mu_j([X])\big |^2
\eeq
we observe that both the absolute value $|\mu_j([X])|$ and $\dim(\ccj)$ are constant on each orbit $\ca_i$. Thus, we may reorganize the sum as:
\begin{eqnarray*}
\dimcc
&=&
\sum_{i=0}^r\sum_{j\in \ca_i}\dimccj\big|\mu_j([X])\big|^2
=\sum_{i=0}^r |\ca_i|
\dim(\cc^\ji) \muji([X])\overline{\muji([X])}
\\&=& 
|\gwkc|\sum_{i=0}^r 
\dim(\cc^\ji) \muji([X])\overline{\muji([X])}
\end{eqnarray*}
The result follows by dividing both sides by $|\gwkc|$ and noting that $\gwkc \simeq \ucc$ by Remark \ref{gwkc-uc}.
\epf
\bl\label{s-geq-0-sgl} 
Let $\cc = \genx$ be a spherical fusion category with commutative Grothendieck ring generated by a simple object $X$. If $\cc$ is $s$-Isaacs, then:
$$
\frac{(\dimcc)^{s+1}}{\oucc\dim(X)}\in \mathbb A.
$$
\el
\bpf
From Equation \eqref{eq-act}, we derive:
\beq\label{eq-act-dimx}
\frac{\dimcc}{|\gwkc|\dim(X)}=\sum_{i=0}^r\dim(\cc^{j_i})\frac{\mu_{\ji}([X])}{\dim(X)}\overline{\mu_{j_i}([X])},
\eeq
Multiplying both sides by ${(\dimcc)^s}/{\dim(X)}$ yields:
\beq\label{eq-act-is-0}
\frac{(\dimcc)^{s+1}}{|\gwkc|\dim(X)}=\sum_{i=0}^r\lambda_0(\mu_\ji, X)\dim(\cc^{\ji})^{s}\overline{\mu_{j_i}([X])}\in \mathbb A.
\eeq
Since $\cc$ is Isaacs, $\lambda_0(\mu_\ji, X)\in \mathbb A$. For $s\geq 0$ we have that $\dim(\cc^{\ji})^{s}\in \mathbb A$ and also $\overline{\mu_{j_i}([X])}$ isn algebraic integer. Then the conclusion follows.
\epf
\subsection{General case for $s\geq 0$.}
We now extend our results to the general setting.
\bt\label{s-geq-0-gf-zcx}
Let $\cc$ be a  spherical fusion category and $X$ be a simple object of $\cc$ such that the fusion subcategory $<X>$ generated by $X$ is $s$- Isaacs fusion category. Then 
$$
\frac{(\dimcc)^{s+1}}{\ndim(Z_\cc(X))\dim(X)}\in \mathbb A.
$$
\et
\bpf
Let $\cd:=\genx$. Applying Lemma \ref{s-geq-0-sgl}  to  $\cd$ gives that:
$$
\dim(X)\mid \frac{(\dimcd)^{s+1}}{\oucd }.
$$
By Equation \eqref{eq-bantay} we have:
\beq\label{eq-bantay-pr}
\frac{(\dimcd)^{s+1}}{\oucd }=\frac{\dimcd}{|U(\cd)|}(\dimcd)^s=\frac{\dimcc}{\ndim(Z_\cc(X))}(\dimcd)^s.
\eeq
Since ${\dimcc}/{\dimcd}\in \mathbb A$ \cite{EGNO15}, we conclude that
$$
\dim(X)\mid \frac{(\dimcc)^{s+1}}{\ndim(Z_\cc(X))}.
$$
\epf
\subsubsection{Proof of Theorem \ref{s-geq-0-gf-uc}}
\begin{proof}
By Theorem \ref{s-geq-0-gf-zcx} we have that $\dim(X)$ divides $\frac{(\dimcc)^{s+1}}{\ndim(\zccx)}$. Using Equation \eqref{ndim} we obtain that $\dim(X)$ divides
\beq\label{eq-bantay-22}
(\dimcc)^{s}\frac{\dimcc}{\ndim(\zccx)}=(\dimcc)^{s}\frac{\fpcc}{n(\zccx)}.
\eeq
The result follows from Theorem \ref{uc-div} since $\wkc$ is real non-negative (RN).
\end{proof}
\subsubsection{Proof of Theorem \ref{ribbon-uc}}
By Remark \ref{br-unitary}, for pseudo-unitary braided fusion categories,  $\wkc$  is real non-negative (RN). Applying Theorem \ref{s-geq-0-gf-uc} with $s = 0$ yields:
$$
\frac{\fpcc}{\oucc \dim(X)}\in \mathbb A.
$$
\bt\label{isaacs-gf-uc-zccx}
Let $\cc$ be a  spherical fusion category and $X$ be a simple object of $\cc$ such that  $\genx$ is Isaacs. Then 
$$\dim(X)\mid \frac{\fp(\cc)}{n(Z_\cc(X))}$$
\et
\bpf
For $s=0$ in Theorem \ref{s-geq-0-gf-zcx} we have 
$$
\frac{\dimcc}{\ndim(Z_\cc(X))\dim(X)}\in \mathbb A.
$$
The conclusion follows from  Equation \eqref{ndim} .
\epf
\bt\label{isaacs-gf-uc}
Let $\cc$ be a spherical fusion category such that $\wkc$ is real non-negative (RN). Let $X$ be a simple object in $\cc$ such that  $\genx$ is Isaacs. Then 
$$\dim(X)\mid \frac{\fpcc}{\oucc}.$$
\et
\bpf
Since $\wkc$ is RN, Theorem \ref{uc-div} implies $\oucc \mid n(\zccx)$. The result then follows from Theorem \ref{isaacs-gf-uc-zccx}.
\epf
\br
Let $\cc$ be a spherical fusion category such that $\genx$ is $s$-Isaacs for some $X\in \irrcc$. 
Let 
$$
C_s:=\frac{(\dimcc)^{2s}\fpcc}{\dim(X)^2\oucc},\;D_s:=\frac{(\dimcc)^{s}\fpcc}{\dim(X)\oucc}.
$$
Note that $$D_s^2=C_s\frac{\fpcc}{\oucc}$$ and this implies that if $C_s\in \mathbb A$ then $D_s\in \mathbb A$ since $\fpcc/\oucc\in \mathbb A$. Thus, if $\genx$ is $s$-Isaac with $s\geq 1/2$ then the divisibility result of Theorem \ref{s-geq-0-gf-uc} is also satsfied.
\er
\section{The modular case}
Let $\cc$ be a premodular fusion category (i.e. braided and spherical). We maintain all previous  notation, particularly $\irr(\cc)=\{X_0, X_1,\dots,X_m\}$ for the simple objects, and $d_i:=\dim(X_i)$ for their quantum dimensions.

By \cite[Example 6.14]{scalg}  there is $\comp$-algebra homomorphism  (the Drinfeld map): $\fq: {\cfcc}\ra \cecc$ given by the following formula:
\beq\label{sh}
\fq(\ch_i)=\sum_{i'=0}^m\frac{s_{ii'}}{d_{i'}}E_{i'}.
\eeq
where $S=(s_{ij})$ is the $S$-matrix of $\cc$.

A premodular fusion category $\cc$ is called {\it modular} if its $S$-matrix is non-degenerated. By \cite[Proposition 3.7]{dgno} the S-matrix is non-degenerate if and only if $\ccpp=\vect$.

In the case of a modular tensor category $\cc$, the Drinfeld map $f_Q:\cfcc\ra\cecc$, as discussed in \cite[Example 6.14]{scalg}, is a an algebra isomorphism. Thus, for any (central) primitive idempotent  $F_j$ of $\cfcc$, we have that $\fq(F_j)$ is a (central) primitive idempotent of $\cecc$. Therefore, we have $\fq(F_j)=E_i$ for some central primitive idempotent $E_i\in \cecc$ corresponding to the simple object $V_i\in \irrcc$.   This yields a bijection $F_j\leftrightarrow V_i$  between the set of (central) primitive idempotents $\{F_j\}$ of $\cfcc$ and $\irrcc$. 

This bijection was also discussed in  \cite[Section 6]{ccc-march}.  It allows us to index the (central) primitive idempotents of $\cfcc$ by the set $\mtc I$ defined above.  With this new indexation, the second item of \cite[Theorem 6.9]{ccc-march} implies that $\dim(\cci)=\dim(V_i)^2$ for any $i \in \mtc I$.

Moreover, under this new notation, applying $\fq^{-1}$ to Equation \eqref{sh} we have that
\beq\label{sh-2}
\ch_i=\sum_{i'=0}^m\frac{s_{ii'}}{d_{i'}}F_{i'}.
\eeq
which implies the crucial evaluation formula:
\beq\label{mf-i}
\mu_{i'}([\ch_i])=\frac{s_{ii'}}{d_{i'}}.
\eeq

\subsection{Proof of Theorem \ref{modular-uc}}
\bpf
Using Equation \eqref{mf-i}  and the fact that $\dim(\cc^i)=d_i^2$ for all $i$,  we transform Equation \eqref{eq-act} as:
\beq\label{eq-act-m}
\frac{\dimcc}{{|\gwkc|}}=\sum_{i=0}^rd_i^2\mu_{X_{i}}([X])\overline{\mu_{X_{i}}([X])}
\eeq
Equation \eqref{mf-i} gives that:
$d_i\mu_{i}([X])=\dim(X)\mu_X([X_i]),$ and Equation \eqref{eq-act-m}  becomes:
\beq\label{eq-act-m2}
\frac{\dimcc}{|\gwkc|\dim(X)^2}=\sum_{i=0}^r\mu_X([{X_{i}}])\overline{\mu_X([{X_{i}}])}
\eeq
which shows that  $\frac{\fpcc}{|\gwkc|\fp(X)^2}\in \mathbb A$ since $\mu_X([{X_{i}}])\in \mathbb A$.
\epf

\subsection{Ito-Michler theorem for modular categories.}
In this section we prove Theorem \ref{im-converse}: 
\bpf
For any simple object $X$, consider $\mathcal{D} = \genx$. Then  $\cd$ is weakly integral and therefore pseudo-unitary. Since $\cd$ is also braided it is Isaacs and therefore Theorem \ref{ribbon-uc} can be applied. 

Thus we have that $\fp(X)^2\mid \frac{\fpcc^2}{\oucc^2}$. Since $\cc$ is weakly integral it follows that $\fp(X)^2\in \mathbb Z_{>0}$ and since $p\nmid  \frac{\fpcc^2}{\oucc^2}$ it follows that $p\nmid \fpx^2$.
\epf
\subsubsection{Proof of Corollary \ref{complete-im}}
It follows from \cite[Proposition 1.3]{b-conj-im} and Theorem \ref{im-converse}.
\subsection{Conjectural Extensions}
In \cite{eno-nec}, the authors proposed the following conjecture:
\begin{conjecture} \label{eno-conj}\cite{eno-nec}
Any spherical fusion category is 1-Isaacs. 
\end{conjecture}  
Using Theorem \ref{s-geq-1/2-gf-uc}, note that  Conjecture \ref{eno-conj} implies the following:
\begin{conjecture}\label{mine}
Let $\cc$ be a spherical fusion category such that $\wkc$ is real non-negative (RN).
For any simple object $X$ of $\cc$, we have:
$$
\frac{(\dimcc)^2\fp(\cc)}{\dimx^2\oucc}\in \mathbb A.
$$
\end{conjecture}
\bibliographystyle{alpha}
\bibliography{24nov}
\ed